\documentclass[12pt]{amsart}
\usepackage{ulem}
\usepackage[colorlinks, bookmarks=true]{hyperref}
\usepackage{latexsym}
\usepackage{amsmath}
\usepackage{amsfonts}
\usepackage{amssymb}
\usepackage{amsthm}
\usepackage{amscd}
\usepackage{color}
\input xy
\xyoption{all} \hfuzz=1.5pt \tolerance=500

\newtheorem{theorem}{Theorem}[section]

\newtheorem{lemma}[theorem]{Lemma}

\newtheorem{proposition}[theorem]{Proposition}

\theoremstyle{definition}
\newtheorem{definition}[theorem]{Definition}
\newtheorem{remark}[theorem]{Remark}

\newtheorem{example}[theorem]{Example}

\newcommand{\ot}{\otimes}

\newcommand{\va}{\varphi}

\newcommand{\id}{{\bf 1}}

\newcommand{\s}{\sigma}
\newcommand{\lm}{\lambda}

\newcommand{\bd}{\begin{document}}
\newcommand{\ed}{\end{document}}
\newcommand{\Htp}{\mathbin{\stackrel{\cdot}{\bigotimes}}}
\newcommand{\Hts}{\mathbin{\stackrel{\cdot}{\bigoplus}}}
\newcommand{\nl}{\nu\in\Lambda}
\newcommand{\xl}{X_\nu;\nl}
\newcommand{\x}{X_\nu}
\newcommand{\la}{\langle}
\newcommand{\ra}{\rangle}
\newcommand{\ch}{\cal H}
\newcommand{\su}{\subseteq}
\newcommand{\B}{\Bbb}
\newcommand{\e}{\varepsilon}
\newcommand{\ov}{\overline}
\newcommand{\vk}{\varkappa}
\newcommand{\vt}{\vartheta}
\newcommand{\al}{\alpha}
\newcommand{\de}{\delta}
\newcommand{\lo}{\Longleftrightarrow}
\newcommand{\T}{{\mathbb T}}
\newcommand{\cK}{{\cal F}}
\newcommand{\sq}{{\square}}
\newcommand{\bb}{{\mathcal B}}
\newcommand{\dd}{{\cal D}}
\newcommand{\cR}{{\cal R}}
\newcommand{\f}{{\cal K}}
\newcommand{\cb}{{\mathcal{CB}}}
\newcommand{\cc}{{\mathcal{CC}}}
\newcommand{\wrr}{\widetilde{\cal R}}
\newcommand{\br}{^\bullet{\cal R}}
\newcommand{\en}{E_\nu}
\newcommand{\el}{{\cal L}}
\newcommand{\lr}{\Longrightarrow}
\newcommand{\bp}{\bigcirc}
\newcommand{\Long}{\Longleftarrow}
\newcommand{\q}{\quad}
\newcommand{\qq}{\qquad}
\newcommand{\cd}{\cdot}
\newcommand{\fur}{f: E\to{\B R}}
\newcommand{\furo}{f_0: E_0\to{\B R}}
\newcommand{\fuc}{f: E\to{\B C}}
\newcommand{\ii}{\infty}
\newcommand{\di}{\diamondsuit}
\newcommand{\bgd}{{\bigtriangledown}}
\newcommand{\bu}{{\bigtriangleup}}
\newcommand{\bc}{{completely bounded}}
\newcommand{\qs}{{quantum space}}
\newcommand{\isc}{{isometric}}
\newcommand{\ism}{{isomorphism}}
\newcommand{\qss}{{quantum spaces}}
\newcommand{\bco}{{completely bounded operator}}
\newcommand{\bcos}{{completely bounded operators}}
\newcommand{\res}{{respectively}}
\newcommand{\tp}{{tensor product}}
\newcommand{\eq}{{equivalent}}
\newcommand{\qtp}{{quantum tensor product}}
\newcommand{\mma}{\mathrel{\mathop{\otimes}\limits_{A}}}
\newcommand{\mmb}{\mathrel{\mathop{\otimes}\limits_{\bb}}}
\newcommand{\mmp}{\mathrel{\mathop{\otimes}\limits_{p}}}
\newcommand{\mmh}{\mathrel{\mathop{\otimes}\limits_{h}}}
\newcommand{\mmop}{\mathrel{\mathop{\otimes}\limits_{4}}}
\newcommand{\mmf}{\mathrel{\mathop{\otimes}\limits_{op}}}
\newcommand{\mmi}{\mathrel{\mathop{\otimes}\limits_{i}}}
\newcommand{\mmm}{\mathrel{\mathop{\otimes}\limits_{\bb-\bb}}}
\newcommand{\mmdd}{\mathrel{\mathop{\otimes}\limits_{\cdot}}}

\newcommand{\mmo}{\mathrel{\mathop{\cdot}\limits_{1}}}
\newcommand{\mmA}{\mathrel{\mathop{\otimes}\limits_{A-A}}}
\newcommand{\mmd}{\mathrel{\mathop{\cdot}\limits_{2}}}
\newcommand{\msp}{\mathrel{\mathop{\otimes}\limits_{sp}}}
\newcommand{\mms}{\stackrel{h}{\otimes}}
\newcommand{\mmt}{\stackrel{4}{\otimes}}
\newcommand{\mmx}{\stackrel{i}{\otimes}}
\newcommand{\mmy}{\stackrel{p}{\otimes}}
\newcommand{\mmz}{\stackrel{sp}{\otimes}}
\newcommand{\gd}{\ddagger}
\newcommand{\od}{\odot}
\newcommand{\mt}{\mapsto}
\newcommand{\mmc}{\mathrel{\mathop{\otimes}\limits_{\sim}}}
\newcommand{\mme}{\stackrel{\sim}{\otimes}}
\newcommand{\spn}{\rm span\,}
\newcommand{\ran}{\mathrm{ran\,}}

\def\query#1{\setlength\marginparwidth{80pt}%
\marginpar{\raggedright\fontsize{10}{10}\selectfont\itshape \hrule\smallskip 
{\textcolor{blue}{#1}}\par\smallskip\hrule}} 
\def\removequeries{\def\query##1{}} 

\newcommand\blfootnote[1]{%
	\begingroup
	\renewcommand\thefootnote{}\footnote{#1}%
	\addtocounter{footnote}{-1}%
	\endgroup
}

\makeatletter
\@namedef{subjclassname@2020}{\textup{2020} Mathematics Subject Classification}
\makeatother

\begin{document}

\numberwithin{equation}{section}

\title[Free and projective generalized multinormed spaces]{Free and projective generalized multinormed spaces}

\author{A.~Ya.~Helemskii}

\address{Faculty of Mechanics and Mathematics, Moscow State (Lomonosov) University, Moscow, Leninskie Gory, Russia 119991} 
\email{helemskii@rambler.ru}

\author{T.~Oikhberg}

\address{ Dept.~of Mathematics, University of Illinois, Urbana IL 61801, USA} 
\email{oikhberg@illinois.edu}



\thanks{The first author was supported by the Russian Foundation for Basic Research (Grant No. 19-01-00447). The second author was supported by the NSF award 1912897.}

% \centerline{{\bf{\Large Free and projective }}} 

% \centerline{{\bf{\Large generalized multinormed spaces}}} 
% \footnote{This paper was written with the support of the Russian Foundation for Basic Research (Grant No. 19-01-00447).} 
% \footnote{
% \blfootnote{ The first author was supported by the Russian Foundation for Basic Research (Grant No. 19-01-00447). The second author was supported by NSF award 1912897.} 

% \bigskip  \centerline{ A.~Ya.~Helemskii and T.~Oikhberg}

\maketitle

\parindent=0pt
\parskip=3pt

\begin{abstract}
The paper investigates free and projective ${\bf L}$-spaces, where ${\bf L}$ is a given normed 
space. These spaces form a far-reaching generalization of known $p$-multinormed spaces; in particular, if 
${\bf L}=L_p(X)$, the ${\bf L}$-spaces can be considered as $p$-multinormed spaces, based on 
arbitrary $\sigma$-finite measure spaces $X$ (for ``canonical'' $p$-multinormed spaces, $X=\mathbb N$ with the counting measure). 
We first describe a ``naturally appearing'' functor, based on paving ${\bf L}$ with contractively complemented finite dimensional subspaces. % (finite dimensionality turns out to be essential here).
This finite dimensionality is essential; it permits us to describe a free ${\bf L}$-space for this functor.
% We deal with two closely related questions: which ${\bf L}$-spaces are free with respect to a functor, taken from a certain class of naturally appearing functors from the category of ${\bf L}$-spaces to the category of sets, and which ${\bf L}$-spaces are projective. 
% Our functors are defined in terms of a fixed family of contractively complemented subspaces of ${\bf L}$. We prove that, for such a functor, free ${\bf L}$-spaces exist if, and only if, all subspaces in question are finite-dimensional (as it is the case with suitably chosen subspaces in $L_p(X)$).
As a corollary,  we obtain a wide variety of projective ${\bf L}$-spaces. For ``nice'' ${\bf L}$ (such as the space of simple $p$-integrable functions on a measure space), we obtain a full description of projective ${\bf L}$-spaces. 
\end{abstract}

\maketitle

% \bigskip  {\bf Keywords:} ${\bf L}$--space, ${\bf L}$--contractive operator, projective ${\bf L}$--space, $\odot$-free ${\bf L}$--space, well composed ${\bf L}$--space. 

% \medskip Mathematics Subject Classification (2020): 46L07, 46M05. 

% \vspace{1cm} 

%  \centerline{{\bf{\large Introduction}}} 
\section*{Introduction}\label{s:intro}

\bigskip
In recent years a new substantial concept appeared in functional analysis. This is the notion of 
$p$-multinormed space~\cite{dal}, a natural extension of multinormed and dual multinormed spaces of 
Dales and Polyakov (where $p$ is 1 or $\ii$) and, on the other hand, of ``operator sequence spaces'' 
of Lambert~\cite{lam} (with $p=2$, and convexity is assumed). These spaces happened to be intimately connected with the 
theory of Banach lattices, that is with the important part of a wide circle of questions concerning 
the concept of positivity. 
More recently, a connection between $2$-multinormed spaces and sectorial operators was discovered and investigated by Kalton, Lorist, and Weis in \cite{KLW}.

After the paper~\cite{dal}, the $p$-multinormed spaces, in their turn, were generalized, in the 
frame-work of the so-called non-coordinate approach, to what can be considered as ``$p$-multinormed 
spaces based on arbitrary measures''~\cite{helizv}. 
Taking ${\bf L}$ is a normed space (for instance, $L_p(X)$ or $L_p^0(X)$, where $X$ is a $\sigma$-finite measure; $L_p^0(X)$ is the span of simple functions in $L_p(X)$), we define an ${\bf L}$-space $E$ as a linear space for which ${\bf L}\ot E$ is equipped with a norm, satisfying certain natural conditions (specified in Section \ref{s:L_spaces}, where relevant definitions and basic facts are collected).
Specializing to $X:=\mathbb{N}$ with the counting measure, we recover $p$-multinormed spaces.

The aim of this paper is to describe projective (= ``homologically best'') ${\bf L}$-spaces.
% in particular for ${\bf L} = L_p(X)$ or $L_p^0(X)$--spaces. 
We consider two versions of this notion. 
The first one is the metric projectivity which is, roughly speaking, a suitable concretization of the general categorical projectivity. The second one is the extreme projectivity, taking into account specific features of functional-analytic structures, based on the presence of norm.
Extremely projective Banach spaces were characterized by Grothendick~\cite{gro}; much later Blecher~\cite{ble} described (in different terms) extremely projective operator spaces. Metric projectivity of Banach spaces appeared, under different names, in old papers of Semadeni~\cite{sem} and Graven~\cite{grav}, and then, for operator spaces, in~\cite{hemsb}. 

% In category theory the main workable approach to study projective objects is based on the notion of free object: if you suggest a reasonable definition of the freeness, then projective objects are exactly retracts of free objects. 
Category theory provides us with a valuable tool for studying projective objects: freeness. 
Our goal is to (i) introduce an appropriate notion of freeness, and (ii) to characterize the free objects. Then projective objects can be described as well, since they are exactly the retracts of free objects. 
This approach has been applied to quite a few categories  in functional analysis, and not only for the metric projectivity, but for the extreme projectivity as well (in the latter case one needs to replace retracts by the so-called near-retracts).
For application to operator spaces and modules see~\cite{hemsb, heruv}; general matricially normed 
spaces are treated in~\cite{hemat}. 

In Section \ref{s:functor} of the present paper we introduce a reasonable notion of a free ${\bf L}$--space. Its definition depends on our choice of a family $L^\nu;\nl$ ($\Lambda$ is an index set) of contractively complemented subspaces which ``paves'' ${\bf L}$.
% such that for every $\nu$ there exists a projection of ${\bf L}$ onto $L^\nu$ of norm one. Such subspaces are called convenient. 
We show that {\it free ${\bf L}$--spaces do exist if, and only if, all our subspaces $L^\nu$ are finite-dimensional} (Theorem \ref{t:criterion for freeness}). If the latter 
condition is fulfilled, we obtain a collection of projective ${\bf L}$--spaces. In particular, this 
is the case when ${\bf L}$ is $L_p(X)$ or $L_p^0(X)$, and every $L^\nu$ is a linear span of 
characteristic functions of several measurable subsets of $X$.  

Addressing freeness, we show that the desirable free ${\bf L}$-spaces are the so-called 
$\oplus_1$-sums of all possible families of some special ${\bf L}$-space, which, in its turn, is 
the $\oplus_1$-sum of all minimal ${\bf L}$-spaces $(L^\nu)^*;\nl$ (Propositions \ref{p:adding 
free} and \ref{p:singleton base}). The words ``minimal ${\bf L}$--space'' mean that we consider in 
${\bf L}\ot(L^\nu)^*$ the injective tensor norm. 
% The basic preparatory construction of the ``if part'' of the indicated criterion is a certain 
One of the importants steps in the proof consists of establishing a bijection, for a given ${\bf L}$-space $E$, between the unit ball of ${\bf L}\ot E$ and a specific set of operators from $(L^\nu)^*$ into $E$. 

In Section \ref{s:projectivity} we turn to the study of projectivity.
% Some applications of this result concern the investigation of projective ${\bf L}$-spaces. 
We call  an ${\bf L}$-space ``well composed'' if it is an $\oplus_1$-sum of a family of ${\bf L}$-spaces, such that each of them is minimal ${\bf L}$--space $L^*$ for some finite-dimensional contractively complemented subspace $L$ of ${\bf L}$ (summands may be different).
Then {\it every retract, respectively near-retract of a well composed ${\bf L}$-space is metrically, respectively extremely projective.}
Moreover, if ${\bf L}$, like in the case of $L^0_p(X)$, is the set-theoretical union of its finite-dimensional contractively complemented subspaces, we obtain a full description of relevant ${\bf L}$-spaces: {\it for such an ${\bf L}$ a given ${\bf L}$-space is metrically, respectively extremely projective if, and only if, it is a retract, respectively near-retract of some well composed ${\bf L}$-space} (Theorem \ref{t:describe_projectivity}).
As a particular case, we obtain (Theorem \ref{t:projective LpX}) a description of metrically and extremely projective $p$-multinormed spaces. 
The same theorem also provides a ``quotient representation'' of $p$-multinormed spaces.
This generalizes and sharpens some of the results, obtained by the second author in \cite{oikPOS}. 
% In the restrictive class of $p$-multinormed spaces, projectivity was investigated by the second author in \cite{oikPOS}. 

\section{${\bf L}$--spaces and ${\bf L}$--bounded operators}\label{s:L_spaces}

First we fix the notation. If $G$ and $E$ are linear spaces, we denote by ${\mathcal{L}}(G,E)$ the set of all linear operators from $G$ to $E$.
If, in addition, $G$ and $E$ are normed, denote by $\bb(G,E)$ the space of all bounded operators between them, endowed with the operator norm, and write $\bb(E)$ instead of $\bb(E,E)$. 
The identity map on 
a set $M$ will be denoted by $\id_M$. The symbol $\ot$ is used for the algebraic tensor product of 
linear spaces and linear operators, and also for elementary tensors. The closed and open unit ball 
in a normed space $E$ is denoted by $B(E)$ and $B^0(E)$, respectively. 

Choose and fix (so far arbitrary) normed space ${\bf L}$, which we shall call {\it base space}. 

\medskip
For an arbitrary base space, ${\bf L}\ot E$ is a left module over the algebra $\bb({\bf L})$ with 
the outer multiplication `` $\cd$ '', well defined by $a\cd(\xi\ot x):=a\xi\ot x$. 

\begin{definition}\label{d:L norm}
% \medskip {\bf Definition 1.1.}
Suppose $E$ is a linear space.
A norm on ${\bf L}\ot E$ is called {\it ${\bf L}$--norm on} $E$, if
%  It is a {\it cross-norm} -- that is, $\|\xi \otimes e\| = \|\xi\| \|e\|$ for any $\xi \in {\bf L}$ and $e \in E$.
 ${\bf L}\otimes E$ is contractive as a left $\bb(L)$-module ${\bf L}\otimes E$: the inequality $\|a\cd u\|\le\|a\|\|u\|$ holds for any $a \in \bb(L)$ and $u \in {\bf L}\ot E$. This will be referred as {\it contractibility property}.
% \end{enumerate}
The space $E$, endowed by an ${\bf  L}$--norm, is called {\it ${\bf L}$--space}. 
\end{definition}

\begin{remark}\label{r:normed}
 Any ${\bf L}$-space $E$ can be equipped with a norm: for $e \in E$, define $\|e\| := \|\xi \otimes e\|$, where $\xi \in {\bf L}$ has norm one. Note that, if $\xi, \eta \in {\bf L}$ have norm one, then, by Hahn-Banach Theorem, there exists a norm one $a \in \bb({\bf L})$ so that $a \xi = \eta$; thus, $\|e\|$ is well-defined.
 With this norm on $E$, the norm on ${\bf L} \otimes E$ becomes a {\it cross-norm} -- that is, $\|\xi \otimes e\| = \|\xi\| \|e\|$ for any $\xi \in {\bf L}$ and $e \in E$.
\end{remark}

% \medskip
 Our principal examples of base spaces ${\bf L}$ are $L_p(X)$, where $X$ is a measure space (always assumed to be $\sigma$-finite), $p\in[1,\ii]$,
 and its subspace $L^0_p(X)$, consisting of simple functions. 
%  as a generalization, to the case of a (possibly) continuous measure, of
As any bounded operator $u : L^0_p(X) \to L^0_p(X)$ extends to a bounded operator $\widetilde{u} : L_p(X) \to L_p(X)$ (of the same norm), any $L_p(X)$-space is also an $L^0_p(X)$-space.
We claim that the converse is also true. Indeed, suppose $E$ is an $L^0_p(X)$-space; denote the corresponding norm on $L^0_p(X) \otimes E$ by $\| \cdot \|_0$. Let $\Lambda$ be the net of all finite $\sigma$-subalgebras of $X$. For $\nu \in \Lambda$, let $L^\nu \subset L^0_p(X)$ be the set of all $\nu$-measurable functions -- that is, of functions $\sum_i \alpha_i \chi_{A_i}$, with scalars $\alpha_i$, and $A_i \in \nu$.
Any such $L^\nu$ is a range of a (contractive) conditional expectation, denoted by $Q^\nu$.
For $u \in L_p(X) \otimes E$, let $\|u\| := \sup_{\nu \in \Lambda} \|Q^\nu \cdot u\|_0$. This quantity is well defined: if $u = \sum_k \xi_k \otimes x_k$, then $\|u\| \leq \sum_k \|\xi_k\| \|x_k\|$. 
% $\|tu\| = |t| \|u\|$, and $\|u+v\| \leq \|u\| + \|v\|$, for any $t \in {\mathbb{R}}$ and $u, v \in L_p(X) \otimes E$.
It is easy to see that $\| \cdot \|$ is a norm which coincides with $\| \cdot \|_0$ on $L^0_p(X) \otimes E$.

Further, $\| \cdot \|$ is a $L_p(X)$-norm -- that is, the inequality $\|a \cdot u\| \leq \|u\|$ holds for any contraction $a \in {\mathcal{B}}(L_p(X))$, and any $u = \sum_{i=k}^n \xi_k \otimes x_k \in L_p(X) \otimes E$.
To establish this, fix  $\varepsilon > 0$, and find $\mu \in \Lambda$ so that $\sum_{k=1}^n \|\xi_k\| \|Q^\mu x_k - x_k\| < \varepsilon$ (this is possible, due to the density of simple functions in $L_p(X)$).
For any $\nu \in \Lambda$, 
$$
Q^\nu (a \cdot u) = \sum_{k=1}^n Q^\nu a \big(Q^\mu\big)^2 \xi_k \otimes x_k + \sum_{k=1}^n Q^\nu a \big( \xi_k - Q^\mu \xi_k \big) \otimes x_k 
$$
(recall that $\big(Q^\mu\big)^2 = Q^\mu$).
As $\| \cdot \|_0$ is an $L^0_p(X)$-norm, we have
$$
\Big\| \sum_{k=1}^n Q^\nu a \big(Q^\mu\big)^2 \xi_k \otimes x_k \Big\|_0 \leq
\big\|Q^\nu a Q^\mu\big\| \Big\| \sum_{k=1}^n Q^\mu \xi_k \otimes x_k \Big\|_0 \leq \|u\| .
$$
From the triangle inequality and Remark \ref{r:normed}, 
$$
\Big\| \sum_{k=1}^n Q^\nu a \big( \xi_k - Q^\mu \xi_k \big) \otimes x_k \Big\|_0 \leq \sum_{k=1}^n \big\|Q^\nu a\big\| \big\| \xi_k - Q^\mu \xi_k \big\| \|x_k\| \leq \varepsilon, 
$$
hence, by the triangle inequality, $\|Q^\nu a \cdot u\|_0 \leq \|u\| + \varepsilon$ for any $\nu \in \Lambda$. Taking the supremum over $\nu$, and recalling that $\varepsilon$ can be arbitrarily small, we conclude $\|a \cdot u\| \leq \|u\|$.

For ${\bf L}$ as above, ${\bf L}$-spaces generalize $p$-multinormed spaces, investigated in e.g.~\cite{dal}. 
Indeed, suppose $X:={\mathbb N}$ with the counting measure. Then $L_p^0(X)$ is $\ell_p^0$ -- the space of finite sequences, with the norm inherited from $\ell_p$. It is easy to see that any $\ell_p^0$-norm on $E$ corresponds to a sequence of cross-norms $\| \cdot \|_n$ on spaces $\ell_p^n \otimes E$, with the property that $\|a \cdot u\|_n \leq \|a\| \|u\|_m$, for any $u \in \ell_p^m \otimes E$, and $u : \ell_p^n \to \ell_p^m$.
The preceding paragraph establishes that any $\ell_p^0$-space is an $\ell_p$-space, and vice versa.

% the notion of an $L_p(X)$-space is obviously equivalent to that of $p$-multinormed space. 

Every linear subspace $F$ of an ${\bf L}$-space $E$ will be considered as an ${\bf L}$-space with 
the ${\bf L}$-norm, induced by the embedding of ${\bf L}\ot F$ into ${\bf L}\ot E$. 

If an operator $\va:G\to E$ between linear spaces is given, we shall use, for the operator $\id_{\bf L}\ot\va:={\bf L}\ot G\to{\bf L}\ot E$, the brief notation $\va_\ii$. Obviously, this is a morphism of left $\bb({\bf L})$-modules. 

% Imitating the definitions, given for various `quantized' structures in~\cite{efr,dal,hemsb}, we give 

The following definition is inspired by various ``quantizations'' (found, for instance, in \cite{efr,dal,hemsb}).

\begin{definition}\label{d:L space}
% \medskip {\bf Definition 1.2.} 
An operator $\va:G\to E$ between ${\bf L}$--spaces is called {\it ${\bf L}$--bounded} if  
$\va_\ii$ is bounded. In a similar way we define the  notions of an {\it ${\bf L}$--contractive 
operator}, {\it ${\bf L}$-coisometric operator} ( = ${\bf L}$-quotient mapping) {\it ${\bf 
L}$-strictly coisometric operator} ( = ${\bf L}$-exact quotient mapping), and so on. In particular, 
our $\va$ is ${\bf L}$-strictly coisometric, respectively ${\bf L}$-coisometric when $\va_\ii$ maps 
$B({\bf L}\ot G)$ onto $B({\bf L}\ot E)$, respectively $B^0({\bf L}\ot G)$ onto $B^0({\bf L}\ot E)$. 
\end{definition}

The operator norm of $\va_\ii$ is denoted by $\|\va\|_{{\bf L}b}$. The set of ${\bf L}$--bounded, 
respectively ${\bf L}$--contractive operators between ${\bf L}$-spaces $G$ and $E$ is denoted by 
$\cb(G,E)$, respectively $\cc(G,E)$. 

\medskip
Let  $P,G,E$ be linear spaces, $\tau:G\to E$, $\va:P\to E$ operators. We recall that an operator  
$\psi:P\to E$ is called a {\it lifting of $\va$ across $\tau$}, if it makes the diagram 
\begin{equation}
\xymatrix@R-10pt@C+15pt{
& G \ar[d]^{\tau}\\
P \ar[ur]^{\psi} \ar[r]^{\va} & E } % \eqno(1)
\label{eq:lifting}
\end{equation}
\noindent commutative. 

\begin{definition}\label{d:projectivity}
% \medskip{\bf Definition 1.3.} 
An ${\bf L}$--space $P$ is called {\it metrically projective}, if for every ${\bf L}$--spaces $G,E$ 
every ${\bf L}$-strictly-coisometric operator $\tau:G\to E$ and every  ${\bf L}$-bounded operator 
$\va:P\to E$ there exists a lifting $\psi:P\to G$ of $\va$ across $\tau$,  such that 
$\|\psi\|_{{\bf L}b}=\|\va\|_{{\bf L}b}$. 

% {\bf Definition 1.4.} 
An ${\bf L}$--space $P$ is called {\it extremely projective}, if for $G,E$ and $\va$ as before, 
every ${\bf L}$-coisometric operator $\tau:G\to E$ and every $\e>0$ there exists a  lifting 
$\psi:P\to G$ of $\va$ across $\tau$, such that $\|\psi\|_{{\bf L}b}<\|\va\|_{{\bf L}b}+\e$. 
\end{definition}

We shall use these definitions in an equivalent form, given by the result below.
% the following equivalent form. Namely, we obviously have 

\begin{proposition}\label{p:characterization of projectivity}
% \medskip {\bf Proposition 1.5.}{\it (i) 
$(i)$ $P$ is metrically projective if, and only if for every ${\bf L}$-strictly-coisometric $\tau:G\to E$ 
and an ${\bf L}$-contractive 
$\va : P \to E$ there exists an ${\bf L}$-contractive lifting of $\va$ across $\tau$. 

$(ii)$ $P$ is extremely projective if, and only if for every ${\bf L}$-coisometric $\tau:G\to E$ and 
an ${\bf L}$-contractive $\va : P \to E$ with $\|\va\|_{{\bf L}b}<1$ there exists an 
${\bf L}$-contractive lifting of $\va$ across $\tau$. $\sq$ 
% )
\end{proposition}

\begin{definition}\label{d:retraction}
% {\bf Definition 1.6.} 
An ${\bf L}$-contractive operator is called {\it retraction}, if it has a right inverse ${\bf 
L}$-contractive operator (which, of course, must be an ${\bf L}$-isometry).
Further, an ${\bf L}$-contractive operator is a {\it  near-retraction}, if, for every $\e>0$, it has a right inverse ${\bf L}$-bounded operator $\rho$ such that  $\|\rho\|_{{\bf L}b}<1+\e$. 
An ${\bf L}$-space $E$ is called {\it retract} ({\it near-retract}) of an ${\bf L}$-space $G$ if there is a retraction (respectively, near retraction) from $G$ onto $E$. 
% Further, we say that $G$ is a {\it near-retract} of $E$ if for any $\e > 0$ there exist a contractive operator from $G$ onto $E$, which has a right inverse $\rho$ with $\|\rho\|_{{\bf L}b}<1+\e$.
\end{definition}

In more geometric terms, $E$ is a retract of $G$ if and only if it is a so-called {\it ${\bf L}$-direct summand of $G$}\footnote{Not to be confused with $L$-summands from the theory of $M$-ideals.}, that is  $E$ is ${\bf 
L}$-isometrically isomorphic to an ${\bf L}$-subspace $F$ of $G$ such that there is an ${\bf 
L}$-contractive projection of $G$ onto $F$.
Following \cite{ble}, we say that $E$ is {\it an almost ${\bf L}$-direct summand of $G$} when for every $\e>0$ there is a subspace $F$ of $E$ such that there exists a projection $Q$ of $G$ onto $F$ with $\|Q\|_{{\bf L}b}<1+\e$ and there exists an ${\bf L}$-topological isomorphism $I:E\to F$ with $\|I\|_{{\bf L}b},\|I^{-1}\|_{{\bf L}b}<1+\e$. 
Equivalently, for any $\e > 0$ there exist a contractive operator from $G$ onto $E$, which has a right inverse $\rho$ with $\|\rho\|_{{\bf L}b}<1+\e$.
Consequently, any near-retract of $G$ must be an almost ${\bf L}$-direct summand.

\begin{proposition}\label{p:retract of projective}
% \medskip {\bf Proposition 1.7.} (i) {\it 
$(i)$ A retract of a metrically projective ${\bf L}$-space is itself metrically projective. 

$(ii)$ A near-retract of an extremely projective  ${\bf L}$-space is itself extremely projective. 
% }
\end{proposition}

\begin{proof}
% \smallskip {\it Proof.} 
We restrict ourselves to the p. (ii), since p. (i) is simpler. Suppose that $P_0$ is extremely 
projective, $P$ is a near-restract of $P_0$, while $\tau$, $\va$, and $G$ are as in Proposition 
\ref{p:characterization of projectivity}(ii). 
% Find $\e \in \big(0, 1/\|\va\|_{{\bf L}b}-1\big)$.
Fix $\lambda \in (\|\va\|_{{\bf L}b}, 1)$, and find an ${\bf L}$-contractive map $\s:P_0\to P$, with a right inverse $\rho$ satisfying $\|\rho\|_{{\bf L}b}<1/\lambda$.
As $\|\va \s\|_{{\bf L}b} \leq \|\va\|_{{\bf L}b} < \lambda$, the extreme projectivity of $P_0$ implies the existence of $\psi_0:P_0\to G$ with  $\tau\psi_0=\va\s$, and $\|\psi_0\|_{{\bf L}b}<\lambda$.
Let $\psi = \psi_0 \rho$; then $\tau \psi = \va \s \rho = \va$ (that is, $\psi$ lifts $\va$). Further, $\|\psi\|_{{\bf L}b} \leq \|\psi_0\|_{{\bf L}b} \|\rho\|_{{\bf L}b} < 1$.
% 
% a near-retraction, while $\tau$, $\va$, and $G$ are as in Proposition \ref{p:characterization of projectivity}(ii). For any $\theta \in(0,1)$ and $\va_0:=\theta\va$ we have $\|\va_0\|_{{\bf L}b}<1$, hence $\|\va_0\s\|_{{\bf L}b}<1$. 
% Therefore there is an ${\bf L}$-contractive $\psi_0:Q\to G$ with  $\tau\psi_0=\va_0\s$. 
% Therefore there exists a $\psi_0:P_0\to G$ with  $\tau\psi_0=\va_0\s$, and $\|\psi_0\|_{{\bf L}b}<1$. Find $\rho:P\to P_0$ with $\|\rho\|_{{\bf L}b} < 1/\|\psi_0\|_{{\bf L}b}$, and $\rho \s = \id_P$.
% $\|\theta^{-1}\psi_0\rho\|_{{\bf L}b}\le\theta^{-1}\|\psi_0\|_{{\bf L}b}\|\rho\|_{{\bf L}b}\le1$. Then, taking 
% 
% Now take $\psi:\theta^{-1}\psi_0\rho$. It is easy to see that $\tau\psi=\va$, and $\|\psi\|_{{\bf L}b} \leq \theta^{-1}$. To complete the proof, recall that $\theta$ can be arbitrarily close to $1$. % $\sq$ 
\end{proof}

\begin{remark}
 The proof above shows that any almost ${\bf L}$-direct summand of an extremely projective  ${\bf L}$-space is extremely projective. 
\end{remark}

\section{The functor $\odot$, and $\odot$-free ${\bf L}$--spaces}\label{s:functor}

Denote by ${\bf L}{\bf Nor_1}$ the category with ${\bf L}$-spaces as objects and ${\bf 
L}$-contractive operators as morphisms. As usually, the category of sets and maps is denoted by 
${\bf Set}$. 

Let $\sq:{\bf L}{\bf Nor_1}\to{\bf Set}$ be (so far) an arbitrary functor. 

% \medskip {\bf Definition 2.1.} 
\begin{definition}\label{d:admissible}
An ${\bf L}$-contractive operator $\tau$ is called {$\sq$-admissible}, if the map $\sq(\tau)$ is 
surjective (i.e. it is a retraction in ${\bf Set}$). 
\end{definition}

\begin{example}\label{e:LNor1}
% \medskip {\bf Example 2.2.} 
Consider the functor $\bigcirc:{\bf L}{\bf Nor_1}\to{\bf Set}$, taking $E$ to   $B({\bf L}\ot E)$, 
and taking the ${\bf L}$-contractive operator $\va : E \to G$ to the restriction of $\va_\ii$ to $B(L \otimes E)$. We see that an ${\bf L}$-contractive operator $\tau$ is 
$\bigcirc$-admissible if, and only if it is ${\bf L}$-strictly coisometric. 
\end{example}

% We pass to our main variety of functors that will be of more use than the previous functor. 
We shall be mostly concerned with a different functor, better suited for describing projective ${\bf L}$-spaces.
% A subspace $L^\nu$ in ${\bf L}$ will be called {\it convenient}, if there exists a projection of norm 1, say $Q^\nu$, from ${\bf L}$ onto $L^\nu$. 
In what follows, we suppose that {\it we are given a family $L^\nu; \nl$, where $\Lambda$ is some index set, of contractively complemented subspaces of ${\bf L}$.} 

\begin{example}\label{e:contractively complemented subspaces}
% \medskip {\bf Example 2.3.} 
We are mainly concerned with ${\bf L} = L^0_p(X)$. Take $\Lambda$ to be the family of all finite sets of pairwise disjoint measurable subsets in $X$; let $L^\nu;\nl$ the span of characteristic functions of sets in $\nu$; $Q^\nu$ is the corresponding conditional expectation. 
% For ${\bf L} = L_p(X)$, we take as $\Lambda$ the set of all finite collections $\nu$ of pairwise disjoint norm one functions $f_1, \ldots, f_N \in L_p(X)$.
% The corresponding projection $Q^\nu$ is defined by $Q^\nu g = \sum_{i=1}^N \big( \int g \widetilde{f}_i ) f_i$, where $\widetilde{f}_i = {\mathrm{sign}} f_i \cdot |f_i|^{p-1}$ for $1 \leq p < \infty$, with a slight modification for $p=\infty$. 
\end{example}

Note that we do not have much latitude in selecting the spaces $L^\nu$. Indeed, any contractive projection on $L^0_p(X)$ extends to that on $L_p(X)$ by continuity. Contractive projections on $L_p(X)$ are known to be ``modified conditional expectations'' (see e.g.~\cite[Theorem 4.3]{Ran}). Further, by \cite{Tzaf}, contractively complemented subspaces of $L_p(X)$ are themselves $L_p$ spaces.

% \medskip
Denote, for brevity, $B_E:=B({\bf L}\ot E), B^0_E:=B^0({\bf L}\ot E), B^\nu_E:=B(L^\nu\ot E), 
 B^{0\nu}_E:=B^0(L^\nu\ot E)$ and consider, for every ${\bf L}$-space $E$, the set 
 ${\bf B}_E:=\textsf{X}\{B^\nu_E;\nl\}$, that 
is the  Cartesian product of the closed unit balls of the normed spaces $B^\nu_E;\nl$. Thus, its 
elements can be represented as `rows' $(...,v_\nu,...)_{\nu \in \Lambda}$, with $v_\nu\in B^\nu_E$. 
Let us introduce the functor 
\[
\odot:{\bf L}{\bf Nor_1}\to{\bf Set},
\]
 taking $E$ to ${\bf B}_E$ and an ${\bf L}$-contractive 
operator $\va: G\to E$ to the map $\odot(\va):\odot(G)\to\odot(E), 
(...,v_\nu,...)\mt(...,\va_\ii(v_\nu),...)$. (This map is well defined since, for every $v\in 
B_G^\nu$, we have $\va_\ii(v)=\va_\ii(Q^\nu\cd v)=Q^\nu\cd\va_\ii(v)\in B_E^\nu)$.

Similarly, we define the space ${\bf B}^0_E:=\textsf{X}\{B^{0\nu}_E;\nl\}$ and introduce the 
functor $\odot^0:{\bf L}{\bf Nor_1}\to{\bf Set}$, taking $E$ to ${\bf B}^0_E$ and an ${\bf 
L}$-contractive operator $\va: G\to E$ to the map $\odot^0(\va):\odot^0(G)\to\odot^0(E)$, 

\begin{proposition}\label{p:coisometric admissible}
% \medskip {\bf Proposition 2.4.} {\it 
Any ${\bf L}$-strictly coisometric operator is $\odot$-admissible, whereas any ${\bf L}$-coisometric 
operator is $\odot^0$-admissible. 
%} 
\end{proposition}

\begin{proof}
% \smallskip {\it Proof.} 
If $\tau:G\to E$ is ${\bf L}$-strictly coisometric, let us take $u\in B^\nu_E$. Then 
$u=\tau_\ii(v)$ for some $v\in B_G$. Hence $\tau_\ii(Q^\nu\cd v)=Q^\nu\cd\tau_\ii v=Q^\nu\cd u=u, 
Q^\nu\cd v\in L^\nu\ot G$ and $\|Q^\nu\cd v\|\le\|Q^\nu\|\|v\|\le1$. The proof in the `${\bf L}$- 
coisometric' case is similar. % $\sq$ 
\end{proof}

In certain cases, the converse of Proposition \ref{p:coisometric admissible} holds.
We say that ${\bf L}$ is {\it properly presented} (relative to the family $L^\nu : \nl$ of complemented subspaces) if for any $N \in \mathbb N$, and $f_1, \ldots, f_N \in {\bf L}$, there exists $\nl$ with $f_1, \ldots, f_N \in L^\nu$.
% if it is the set-theoretical union  of a family of its contractively complemented {\it finite-dimensional} subspaces.
Proper presentation occurs, for instance, in the setting of Example \ref{e:contractively complemented subspaces}. For $f_1, \ldots, f_N \in {\bf L} = L^0_p(X)$, let $\nu$ be the (finite) $\sigma$-algebra of subsets of $X$, generated by $f_1, \ldots, f_N$. Clearly $f_1, \ldots, f_N \in L^\nu$.

\begin{proposition}\label{p:odot criterion}
% \medskip  {\bf Proposition 2.5.} {\it 
Suppose ${\bf L}$ is properly presented relative to the family $L^\nu : \nl$, and let $\odot$ and $\odot^0$ be the functors arising from this family. Suppose, further, that $G$ and $E$ are ${\bf L}$-spaces, and $\phi : G \to E$ is a linear operator. Then:
\begin{enumerate}
 \item 
 $\|\phi\|_{{\bf L}b} \leq 1$ if and only if $\odot \phi$ $($equivalently, $\odot^0 \phi)$ is a well-defined map on ${\bf Set}$.
 \item 
 $\phi$ is $\odot-$, respectively $\odot^0$--admissible if, and only if, it is  ${\bf L}$-strictly coisometric, respectively ${\bf L}$-coisometric. 
 $\sq$
\end{enumerate}
  %} 
 \end{proposition}

%\medskip
%{\bf Remark 2.6.} We do not know whether the converse is true. Of course, if we replace, as our 
%base state, ${\bf L}$ by ${\bf L}^0$, then our proposition turns out to be a criterion of the 
%admissibility. But this does not affect our future assertions concerning ${\bf L}$-spaces. 

\medskip
The following definition is a particular case of a well known general-categorical definition (see, 
e.,g., Definition 6 in\cite{hemat}). In what follows, $\sq:{\bf L}{\bf Nor_1}\to{\bf Set}$ is an 
arbitrary functor, $M$ a set. 

\smallskip

\begin{definition}\label{d:free space}
% \medskip {\bf Definition 2.6.} 
An ${\bf L}$-space ${\bf F}^\sq(M)$ is called a {\it $\sq$-free space with the base $M$}, if,  for 
every ${\bf L}$-space $E$ there exists a bijection 
$$
{\bf I}_E: {\bf Set}(M,\sq E)\to\cc({\bf F}^\sq(M),E) %,\eqno(2)
$$
between the respective sets of morphisms, and these bijections have the so-called {\it natural 
property.} The latter means that for every $\va\in\cc(G,E)$ we have the commutative diagram 
\begin{equation}
% $$
\xymatrix@C+20pt{ {\bf Set}(M,\sq G) \ar[r]^{ {\bf I}_G} \ar[d]_{(\sq\va)^*}
& \cc({\bf F}^\sq(M),G) \ar[d]^{\va^*} \\
{{\bf Set}(M,\sq E)} \ar[r]^{{{\bf I}_E}} & \cc({\bf F}^\sq(M),E).)}
% \eqno(2)
% $$
\label{eq:natural property}
\end{equation}
where $\va^*$ acts by composition as $\psi\mt\va \circ \psi$ and $(\sq\va)^*$ as $\rho\mt(\sq\va) \circ \rho$. 
\end{definition}

A simple ``diagram-chasing'' argument shows that the free space of $M$ is unique.

Henceforth, we shall write ${\bf F}(M)$ for ${\bf F}^\sq(M)$, if there is no confusion about the functor $\sq$.

%\noindent \fbox{\begin{minipage}{\textwidth} T.O.: On horizontal arrows, I removed 
%squares from ${\bf I}^\sq_G$ and ${\bf I}_E^\sq$, for the sake of consistency with the preceding 
%definition, and subsequent usage.
%\end{minipage}}

\medskip
The following proposition is well known in the general categorical context of adjoint functors 
~\cite{mcl}. 

\begin{proposition}\label{p:quotient onto base}
% \medskip {\bf Proposition 2.7.} {\it 
Suppose that for a given ${\bf L}$-space $E$ there exists a $\sq$-free space ${\bf F}(\sq E)$ with 
the base $\sq E$. Then there exist a $\sq$-admissible operator $\pi: {\bf F}(\sq E)\to E$. 
% } 
\end{proposition}

\begin{proof}
% \smallskip {\ Proof.} 
% Indeed, for $\pi:={\bf I}_E(\id_{\sq E})$ the natural property gives 
% \begin{align*}
% \sq\pi \circ {\bf I}_{{\bf F}(\sq E)}^{-1}(\id_{{\bf F}(\sq E)})
% &
% =
% {\bf I}_E^{-1}{\bf I}_E \big( \sq\pi \circ {\bf I}_{{\bf F}(\sq E)}^{-1}(\id_{{\bf F}(\sq E)}) \big)
% \\
% &
% =
% {\bf I}_E^{-1}\big[\pi \circ {\bf I}_{{\bf F}(\sq E)}{\bf I}_{{\bf F}(\sq E)}^{-1}(\id_{{\bf F}(\sq E)}) \big]=
% {\bf I}^{-1}_E(\pi)=\id_{\sq E}. % \qedhere % \q \sq
% \end{align*}
%
To simplify the notation, denote ${\bf F}(\sq E)$ by $G$. We show that $\pi:={\bf I}_E(\id_{\sq E})$ is admissible, by establishing that $\sq \pi \circ \psi = \id_{\sq E}$, where $\psi = {\bf I}_G^{-1}(\id_G)$. To this end, observe that, due to the natural property \eqref{eq:natural property},
$$
{\bf I}_E \big( \sq \pi \circ \psi \big) = \pi \circ {\bf I}_G \big( {\bf I}_G^{-1}(\id_G) \big) = \pi \circ \id_G = \pi ,
$$
hence
$$
\sq \pi \circ \psi = 
{\bf I}_E {\bf I}_E^{-1} \big( \sq \pi \circ \psi \big) = {\bf I}_E^{-1} \pi = \id_{\sq E} . \qedhere
$$
\end{proof}

The main result of the rest of this section -- Theorem \ref{t:criterion for freeness} -- determines 
the necessary and sufficient condition for the existence of $\odot$-free ${\bf L}$-spaces. 

% \medskip
% When, for functor $\odot$, there exist respective $\odot$-free ${\bf L}$-spaces? The answer is our future Theorem \ref{t:criterion for freeness}, which we begin to prepare. 
% \medskip
First we establish that only the case of finite dimensional spaces $L^\nu$ is interesting.

%\noindent \fbox{\begin{minipage}{\textwidth} T.O.: This is the former Proposition 
%2.13. I moved it here, in order to rule out the infinite dimensional case, and focus on the finite 
%dimensional one in the sequel. 
%\end{minipage}}

 \begin{proposition}\label{p:must be fin dim}
% \medskip {\bf Proposition 2.13.} {\it 
Suppose that, for some set $M$, there exists a $\odot$--free ${\bf L}$--space with the base $M$. 
Then all subspaces $L^\nu$ are finite-dimensional. 
% }
\end{proposition}

\begin{proof}
% \smallskip {\it Proof.} 
Let ${\bf F}$ be our $\odot$--free ${\bf L}$--space. Taking, in the capacity of $E$, the same ${\bf 
F}$, we obtain a bijection ${\bf I}_{\bf F}:{\bf Set}(M,\odot{\bf F})\to\cc({\bf F},{\bf F})$. 
% Since the operator $\id_{\bf F}$ is, of course, ${\bf L}$--contractive, we can speak about the 
%function $f:M\to{\bf B}_{\bf F}$, which is ${\bf I}_{\bf F}^{-1}(\id_{\bf F})$. 
Let $f = {\bf I}_{\bf F}^{-1}(\id_{\bf F})$. For any  $t\in M$, $f(t)$ is a 
``row,'' say $(...,u_\nu,...);\nl$, where $u_\nu\in B^\nu_{\bf F}$ for all $\nu$. 

% From now on we choose an arbitrary $\nl$ and fix it.
Fix $\nl$, and represent $u_\nu \in L^\nu\ot{\bf F}$ as $\sum_{k=1}^n\xi_k\ot x_k$ 
for some $\xi_k\in L^\nu, x_k\in{\bf F};  k=1,...,n$. Take an arbitrary non-zero $\eta\in L^\nu$ 
and $x\in{\bf F}$ such that $\|\eta\ot x\|\le1$. Consider an arbitrary function $g:M\to{\bf B}_{\bf 
F}$ where $g(t)\in{\bf B}_{\bf F}$ is a `row' with $\eta\ot x$ on the `$\nu$-th' place. 
Then $\va := {\bf I}_{\bf F}(g)$ is an ${\bf L}$--contractive operator, acting on 
${\bf F}$.
%say $\va$. 
Therefore, by the ``natural property'', we have 
\[
{\bf I}_{\bf F}(\odot\va \circ f)=\va{\bf I}_{\bf F}(f)=\va{\bf I}_{\bf F}{\bf I}_{\bf F}^{-1}(\id_{\bf F})=
\va\id_{\bf F}=\va={\bf I}_{\bf F}(g) ,
\]
hence $\odot\va \circ f=g$. % $(\odot\va)(f)=g$.
This, in particular, means that the `row' $[\odot\va \circ f](t)$ 
coincides with the `row' $g(t)$. Looking at the `$\nu$-th' place in these `rows', we see that 
$\va_\ii(u_\nu)=\eta\ot x$, that is $\eta\ot x=\sum_{k=1}^n\xi_k\ot\va(x_k)$. This % obviously 
implies that $\eta\in \spn\{\xi_k :1 \leq k \leq n\}$.
As $\eta \in L^\nu$ is arbitrary, the space  $L^\nu$ coincides with the latter span.
% But $\eta$ and $x$ were arbitrarily chosen. Consequently,   the space $L^\nu$ coincides with the latter span, and hence it is finite-dimensional. % $\sq$ 
 \end{proof}

 To deal with the functor $\odot$, we shall use a certain family of simpler functors. Fix, 
for a time, $\nl$ and consider the functor $\odot_\nu$, taking $E$ to $B^\nu_E$ and taking an ${\bf 
L}$-contractive operator $\va: G\to E$ to the respective restriction 
$\odot(\va):\odot(G)\to\odot(E)$ of the operator $\va_\ii$. Thus, for every $\va\in\cc(E,F)$ and 
${\bf u}=(...,u_\nu,...)\in{\bf B}_G$ we have 
\begin{equation}
(\odot\va({\bf u}))_\nu=\odot_\nu\va(u_\nu). % \eqno(3) 
\label{eq:odot nu}
\end{equation}

Fix a projection $Q^\nu = Q:{\bf L}\to L^\nu$ of norm 1. Equip the dual space $(L^\nu)^*$ with the so-called minimal ${\bf L}$-norm. That is, we supply ${\bf L}\ot(L^\nu)^*$ with the norm of the injective tensor product of respective normed spaces (see, e.g.,~\cite{rya}). For an ${\bf L}$-space $E$ we introduce the operator 
\[ 
{\bf I}^\nu_E:{\bf L}\ot E\to{\mathcal L}((L^\nu)^*,E), 
\]
assigning to $u\in{\bf L}\ot E$ the operator, 
taking $g \in (L^\nu)^*$ to $(g\ot\id_E)(Q\cd u)$. 
%(Cf. a prototype of this construction, presented 
%in the case of $p$-multinormed spaces in~\cite{oif}). 

In Propositions \ref{p:norm est}, \ref{p:isometric isomorphism}, and \ref{p:singleton base} below 
we suppose that {\it the space $L^\nu$ is finite-dimensional.} In this case we take an Auerbach 
basis in $L^\nu$, that is a basis $e_k^\nu = e_k; k=1,...,n_\nu$ in that space together with a basis $e_k^{\nu*} = e_k^*$ in 
$(L^\nu)^*$ such that $e^*_k(e_l)=\de_{kl}$, and $\|e_k\| = 1 = \|e_k^*\|$ for any $k$ (see, e.g., \cite[p. 12]{ltz}). We fix a distinguished element 
\begin{equation}
w^\nu = w:=\sum_ke_k\ot e_k^*\in{\bf L}\ot(L^\nu)^* .
\label{eq:define w}
\end{equation}
As before, $Q^\nu$ (or $Q$ for short) is a fixed projection from ${\bf L}$ onto $L^\nu$.

\begin{proposition}\label{p:norm est}
% \medskip {\bf Proposition 2.8}. {\it 
For every $u\in{\bf L}\ot E$ the operator ${\bf I}^\nu_E(u)$ is ${\bf L}$--bounded, and we have 
$\|{\bf I}^\nu_E(u)\|_{{\bf L}b}= \|Q\cd u\|$. 
% } 
\end{proposition}

\begin{proof}
% \smallskip {\it Proof}. 
Our task is to show that, for any  $u\in{\bf L}\ot E$, we have $\|{\bf I}_E^\nu(u)_\ii\|=\|Q\cd u\|$. %Let us 
Define the linear map $J: {\bf L}\ot(L^\nu)^*\to\bb(L^\nu,{\bf L})$ by letting $J(\xi \ot g)$ ($\xi \in {\bf L}, g \in (L^\nu)^*$ be the operator $\eta\mt g(\eta)\xi$ ($\eta \in L^\nu$).
As ${\bf L}\ot(L^\nu)^*$ is endowed with the injective tensor product norm, $J$ is an isometry (see, e.g., \cite{rya}).
% Recall that norm on ${\bf L}\ot(L^\nu)^*$ is the norm of the injective tensor product of respective normed spaces. In this case we know that the canonical operator $J: {\bf L}\ot(L^\nu)^*\to\bb(L^\nu,{\bf L})$, well defined by taking an elementary tensor $\xi\ot g$ to the  operator, acting as  $\eta\mt g(\eta)\xi$, is an isometry (see, e.g., \cite{rya}). But obviously we 
Further, $J(w)(e_l)=e_l$ for all $l=1,...,n$. Therefore $J(w)$ is just the canonical embedding of $L^\nu$ into ${\bf L}$, and hence $\|w\|=1$. 

For $u:=\xi\ot x$ we have ${\bf I}^\nu_E(u)(e_k^*)=e_k^*(Q\xi)x$. By linearity, ${\bf I}^\nu_E(u)_\ii(w)=\sum_ke_k\ot e_k^*(Q\xi)x$. The equality $\eta = \sum_k e_k^*(\eta) e_k$ is true for any $\eta \in L^\nu$,
and therefore,
$${\bf I}^\nu_E(u)_\ii(w)=\sum_ke_k\ot e_k^*(Q\xi)x = Q \xi \ot x = Q \cd u $$
holds for any elementary tensor product $u$.
By linearity, ${\bf I}^\nu_E(u)_\ii(w)=  Q \cd u $ for any $u\in{\bf L}\ot E$.
% whereas $Q\cd u=Q\xi\ot x$.
% But $Q\xi$ is a linear combination of elements $e_l$. Clearly, for $Q\xi:=e_l$ both expressions transform to $e_l\ot x$. Therefore we have ${\bf I}^\nu_E(u)_\ii(w)=Q\cd u$ for all elementary tensors in ${\bf L}\ot E$, hence for all $u\in{\bf L}\ot E$.
To show $\|{\bf I}^\nu(u)_\ii\|\ge\|Q\cd u\|$, recall that $\|w\| \leq 1$. 

To prove the converse, observe first that $J(v) \cd (Q\cd u)={\bf I}^\nu_E(u)_\ii(v)$ holds for any $u \in{\bf L}\ot E$ and $v \in {\bf L} \otimes (L_\nu)^*$.
By linearity, it suffices to consider elementary tensors $u = \xi\ot x$ and $v = \eta \ot f$.
It is easy to check that then, $J(v) \cd (Q\cd u)=f(Q\xi)\eta\ot x$, whereas ${\bf I}^\nu_E(u)_\ii(v) = \eta\ot f(Q\xi)x$.
% holds for any $v\in{\bf L}\ot(L^\nu)^*;\|v\|\le1$. Indeed, we only need to establish this identity for elementary tensors $\xi\ot x\in{\bf L}\ot E$ and $\eta\ot f\in{\bf L}\ot(L^\nu)^*$
% take an arbitrary $v\in{\bf L}\ot(L^\nu)^*;\|v\|\le1$ and note that  we have 
% \[ (J(v)Q)\cd u={\bf I}^\nu_E(u)_\ii(v). \]
% (It can be easily checked on the elementary tensors of the form $\xi\ot x\in{\bf L}\ot E$ and 
% $\eta\ot f\in{\bf L}\ot(L^\nu)^*$: indeed, the left part transforms to  $f(Q\xi)\eta\ot x$ whereas  the right part to $\eta\ot f(Q\xi)x$.) 
Therefore, taking into account the contractibility property 
and our choice of norm on ${\bf L}\ot(L^\nu)^*$, we have 
\[ 
 \|{\bf I}^\nu_E(u)_\ii(v)\|\le\|J(v)\cd(Q\cd u)\|\le\|J(v)\|\|Q\cd u\|=\|Q\cd u\|\|v\|
\] 
Taking the supremum over all $v$ of norm not exceeding $1$, we conclude that $\|{\bf I}^\nu_E(u)_\ii\|\le\|Q\cd u\|$. \q % $\sq$ 
\end{proof}

\begin{proposition}\label{p:isometric isomorphism}
% \medskip {\bf Proposition 2.9.} {\it 
The operator ${\bf I}^\nu_E$, being restricted to $L^\nu\ot E$, is an  ${\bf L}$-isometric 
isomorphism between the latter space and $(\cb((L^\nu)^*,E))_{{\bf L}b}$. 
% } 
\end{proposition}

\begin{proof}
% \smallskip {\it Proof.} 
Since elements of $L^\nu\ot E$ are exactly those $u\in{\bf L}\ot E$ with $u=Q\cd u$, Proposition 
\ref{p:norm est} implies that our restriction of ${\bf I}^\nu_E$ is ${\bf L}$-isometric.
% To prove that it is a bijection, we shall display its inverse. Indeed, it is easy to check that the operator  $\bb((L^\nu)^*,E))\to L^\nu\ot E$, taking $\va$ to $\va_\ii(w)$ is what we need.% $\sq$  
It is easy to check that the inverse map is given by $\cb((L^\nu)^*,E))\to L^\nu\ot E : \va \mapsto \va_\ii(w)$; this establishes the bijectivity.
\end{proof}

From now on let us agree to {\it identify every map from a singleton to some set with its image in 
this set.} %We pass to the case of an arbitrary family of contractively complemented subspaces $L^\nu;\nl$. 
 
 \begin{proposition}\label{p:singleton base}
% \medskip {\bf Proposition 2.10}. {\it 
The space $(L^\nu)^*$, supplied with the minimal ${\bf L}$-norm, is the $\odot_\nu$-free ${\bf L}$-space whose base is a singleton.
% There exists a $\odot_\nu$-free ${\bf L}$-space with the base a singleton. It is the space $(L^\nu)^*$, supplied with the minimal ${\bf L}$-norm. 
% } 
\end{proposition}

\begin{proof}
% Retain the notation ${\bf I}^\nu_E$ for the restriction of the operator from the previous proposition to $B^\nu_E$. Obviously, it is a bijection between $B^\nu_E$ and the set $\cc((L^\nu)^*,E)$.
Proposition \ref{p:isometric isomorphism} shows that ${\bf I}^\nu_E$ implements a bijection between $B^\nu_E$ and the set $\cc((L^\nu)^*,E)$.
As to the natural property, it is easy to check that for every $u\in B^\nu_G$ 
and $\va\in\cc(G,E)$ both maps $\va^*({\bf I}^\nu_G(u))$ and ${\bf I}^\nu_E(\odot_\nu \va)^*(u)$ (cf. diagram \eqref{eq:natural property}) take $f\in(L^\nu)^*$ to $(f\ot\va)u$. 
\end{proof}

For the rest of this paper, we shall use the following construction. For a family $(E_i)_{i \in \Delta}$ of ${\bf L}$-spaces, 
% where $i$ runs an arbitrary index set, say $\Delta$, 
consider their {\it algebraic} sum $\oplus_i E_i$ and make it an ${\bf L}$-space in the following 
way: we identify ${\bf L}\ot(\oplus_i E_i)$ with $\oplus_i({\bf L}\ot E_i)$ and consider in the 
latter space the norm of the (non-completed) $l_1$-sum of its direct summands (that is, for an 
element $u\in\oplus_i({\bf L}\ot E_i); u=(...,u_i,...);u_i\in {\bf L}\ot E_i$ we set 
$\|u\|_1:=\sum_i\|u_i\|)$. Obviously, we obtain an ${\bf L}$-space, called the {\it $\oplus_1$-sum} 
of our spaces and denoted by $(\oplus_i E_i)_1$. Together with the latter, we consider the family 
of operators ${\bf in}_i:E_i\to(\oplus_i E_i)_1:x\mt(...0,0,x_i,0,0,..);x_i=x$; $i\in\Delta$. These 
are, of course, ${\bf L}$-isometric and hence ${\bf L}$-contractive. 

Now suppose that for some ${\bf L}$-space $E$ we are given a family of ${\bf L}$-contractive 
operators $\va_i:E_i\to E;i\in\Delta$. Consider the operator $\oplus_i\va_i:(\oplus_i E_i)_1\to 
E:x=(...,x_i,...)\mt\sum_i\va_i(x_i)$, well defined by $(\oplus_i\va_i){\bf{in}}_i=\va_i$. It is 
easy to see that $\oplus_i\va_i$ is ${\bf L}$-contractive, and the map $\{\va_i; 
i\in\Delta\}\mt\oplus_i\va_i$ is a bijection between the set of all families 
$\{\va_i\in\cc(E_i,E)\}$ and the set $\cc(\oplus_i E_i,E)$. 

In the categorical language our construction means that the space $(\oplus_i E_i)_1$, together 
with the family ${\bf in}_i$, is the coproduct of objects $E_i$ in ${\bf L}{\bf Nor_1}$.

It is easy to see that for given $\va_i$ and a ${\bf L}$-contractive operator $\psi:E\to F$, where 
$F$ is another ${\bf L}$-space, we have 
\begin{equation}
 \psi(\oplus_i\va_i)=\oplus_i(\psi\va_i)  %\eqno(4)
 \label{eq:oplus operators} 
\end{equation}
Moreover, we obviously have 
\begin{equation}
\|\oplus_i\va_i\|_{{\bf L}b}=\sup_i\|\va_i\|_{{\bf L}b} % \eqno(5) 
\label{eq:oplus norms} 
\end{equation}

%\noindent \fbox{\begin{minipage}{\textwidth} T.O.: Shall we stick to finite 
%dimensional spaces $L^\nu$ until we arrive at Theorem \ref{t:criterion for freeness}? Or if the 
%goal is to develop some general theory, perhaps we should use notation different from $\odot_\nu$, 
%whose meaning is already fixed? 
%\end{minipage}}

From now on, we shall denote the $\odot_\nu$-free space of a singleton, constructed in Proposition \ref{p:singleton base}, by ${\bf F}^\nu(\star)$

\begin{proposition}\label{p:adding free}
%Suppose that for every $\nl$ there exists a $\odot_\nu$-free ${\bf L}$-space with a singleton, say 
%$\{\star\}$, as its base. Then, if ${\bf F}^\nu(\star)$ is the mentioned space, 
%
%\smallskip
$($i$)$
${\bf F}(\star):=(\oplus_\nu{\bf F}^\nu(\star))_1$ is a $\odot$-free ${\bf L}$-space with the base $\{\star\}$.
% there exists a $\odot$-free ${\bf L}$-space with the base $\{\star\}$,   and it is ${\bf F}(\star):=(\oplus_\nu{\bf F}^\nu(\star))_1$.

\smallskip 

$($ii$)$
For any set $M$, ${\bf F}(M):=(\oplus_{t\in M}{\bf F}(t))_1$ is a $\odot$-free  ${\bf L}$-space with $M$ as its base. Here, ${\bf F}(t)$ is ${\bf F}(\star)$ for $t:=\star$. 
% For every set $M$ there exists a $\odot$-free  ${\bf L}$-space with $M$ as its base, and it is ${\bf F}(M):=(\oplus_{t\in M}{\bf F}(t))_1$, where ${\bf F}(t)$ is ${\bf F}(\star)$ for $t:=\star$. 
% } 
\end{proposition}

\medskip In the categorical language p. (ii) means that ${\bf F}(M)$ is the coproduct of 
$card(M)$ copies of ${\bf F}(\star)$. 

\begin{proof}
% \medskip {\it Proof}. 
(i) Identify ${\bf Set}(\{\star\},B^\nu_E)$ with $B^\nu_E$.
Proposition \ref{p:norm est} shows that, for any ${\bf L}$-space $E$, the map ${\bf I}^\nu_E:{\bf Set}(\{\star\},B^\nu_E)\to\cc({\bf F}^\nu(\star),E)$ is bijective.

% Let, for an arbitrary ${\bf L}$-space $E$, ${\bf I}^\nu_E:{\bf Set}(\{\star\},B^\nu_E)\to\cc({\bf F}^\nu(\star),E)$, that is a map from $B^\nu_E$ onto $\cc({\bf F}^\nu(\star),E)$, be the relevant bijection.  

For ${\bf u}=(...,u_\nu,...)\in{\bf B}_E=\odot(E)$ define ${\bf I}^\star_E({\bf 
u}):=\oplus_\nu({\bf I}^\nu_E)(u_\nu)\in\cc({\bf F}(\star),E)$. 
% , where the summantion runs over all $\nl$. 
Since ${\bf I}_E^\nu$ is a bijection for every $\nl$, the same is obviously true for 
${\bf I}^\star_E$. % Finally, the latter bijections have the natural property. Indeed, for 
To establish the natural property, fix ${\bf L}$-spaces $E$ and $G$, and $\va\in\cc(G,E)$.
Pick ${\bf u}=(...,u_\nu,...)\in{\bf B}_G$. The natural property of bijections ${\bf I}^\nu_E$, together with  \eqref{eq:odot nu} and \eqref{eq:oplus operators}, yields
\begin{align*}
\va{\bf I}^\star_G({\bf u})
&
=
\va[\oplus_\nu{\bf I}^\nu_G(u_\nu)]=\oplus_\nu[\va({\bf I}^\nu_G(u_\nu))]
\\
&
=
\oplus_\nu{\bf I}^\nu_E[(\odot_\nu\va)(u_\nu)]=\oplus_\nu{\bf I}^\nu_E[(\odot\va)(u)]_\nu=
{\bf I}_E^\star(\odot\va({\bf u})).
\end{align*}

(ii) Let $E$ be as before. Suppose a function $f:M\to\odot(E)$, or, what is the same, $f:M\to{\bf B}_E$, is given. Set ${\bf I}_E(f):=\oplus_{t\in M} {\bf I}_E^t(f(t))\in\cc({\bf F}(M),E)$, where ${\bf I}_E^t$ ($t \in M$) are copies of ${\bf I}_E^\star$.
The map ${\bf I}_E:{\bf Set}(M,\odot(E)) \to\cc({\bf F}(M),E)$ is bijective, since all ``summands'' ${\bf I}_E^t$ are. The relevant natural property follows, for $\va\in\cc(G,E)$ and $f:M\to{\bf B}_G$, from the equalities 
\[
\va{\bf I}_G(f)=\va[\oplus_t{\bf I}_G^t(f(t))]=\oplus_t[\va{\bf I}_G^t(f(t))]=
\oplus_t{\bf I}_E^t(\odot\va)(f(t))={\bf I}_E(\odot\va(f)). %\q \sq
\qedhere
\]
\end{proof}

\begin{remark} \label{r:general free}
Actually, Proposition \ref{p:adding free} is but a particular case of a certain general-categorical assertion. 
Namely, let ${\mathcal K}$ be some category with coproduct. We shall denote the coproduct of a family 
$Y_i;i\in\Delta$ of objects in ${\mathcal K}$ by $\oplus_iY_i$. Further, consider a family of functors 
$\sq_\mu:{\mathcal K}\to{\bf Set}$, where $\mu$ runs some index set $\Lambda'$. Introduce the ``combined'' functor $\sq:{\mathcal K}\to{\bf Set}$, taking $Y$ to the Cartesian product of sets 
$\sq_\mu(Y);\mu\in\Lambda'$. Further, $\sq$ takes a morphism $\va:Y\to Z$ to the map, transforming a ``row'' $(...,u_\mu,...);u_\mu\in\sq_\mu(Y)$ to $(...,\sq_\mu\va(u_\mu),...)$.
Now suppose that for every  $\mu$ there exists a $\sq_\mu$-free object in ${\mathcal K}$ with a singleton, say $\{\star\}$, as its base. Then, if ${\bf F}^\nu(\star)$ is the space constructed above, then: 

(i) ${\bf F}(\star):=(\oplus_\mu{\bf F}^\mu(\star))$ is a $\sq$-free  object in ${\mathcal K}$ with the base $\{\star\}$.
% there exists a $\sq$-free  object in ${\mathcal K}$ with the base $\{\star\}$,   and it is ${\bf F}(\star):=(\oplus_\mu{\bf F}^\mu(\star))$.

(ii) For every set $M$ there exists a $\sq$-free  object in ${\mathcal K}$ with $M$ as its base, and it 
is ${\bf F}(M):=(\oplus_{t\in M}{\bf F}(t))$, where ${\bf F}(t)$ is ${\bf F}(\star)$ for 
$t:=\star$. 

The proof proceeds as in Proposition \ref{p:adding free}, with obvious modifications.
\end{remark}

% {\bf Remark 2.12.} 
\begin{remark}\label{r:a kind of freeness}
Consider the category ${\bf L}{\bf Nor}$ of ${\bf L}$-spaces and all ${\bf L}$-bounded operators. Define the functor $\overline{\odot} : {\bf L}{\bf Nor} \to {\bf Set}$, taking $E$ to $L^\nu\ot E$.
Modifying the proof of Proposition \ref{p:singleton base}, one shows that $(L^\nu)^*$ is a free ${\bf L}$-space with respect to $\overline{\odot}$. However, this observation does not produce a full analogue of Proposition \ref{p:adding free}(ii): 
%\query{Should parts (ii) and (iii) of Theorem \ref{p:singleton base} be parts (i) and (ii) of 
%Theorem \ref{p:adding free}?} 
% The obstacle is that the 
the category ${\bf L}{\bf Nor}$, unlike ${\bf L}{\bf Nor_1}$, has no coproducts of infinite families of ${\bf L}$-spaces. 
\end{remark}
 
Combining Propositions \ref{p:singleton base}, \ref{p:adding free} and \ref{p:must be fin dim}, we 
immediately obtain 

\begin{theorem}\label{t:criterion for freeness}
% \medskip {\bf Theorem 2.14.} 
For every set $M$ there exists a $\odot$-free ${\bf L}$--space with the base 
$M$ if, and only if, all spaces $L^\nu$ are finite-dimensional. %Moreover,  
\end{theorem}

\begin{remark}
% \medskip {\bf Remark 2.15.} 
In particular, we see that for the functor $\bigcirc$ from Example \ref{e:LNor1} $\bigcirc$-free 
${\bf L}$--spaces exist only when ${\bf L}$ is finite-dimensional. In this setting, our family of subspaces $(L^\nu)$ consists of only one space, namely ${\bf L}$ itself. 
\end{remark}

\medskip
For $f:M\to{\bf B}_E:t\mt {\bf u}(t):=(...,u_\nu(t),...)$ we define 
$|||f|||:=\sup_t\sup_\nu\| u_\nu(t)\|$ and, for $\lm\in[0,1]$, 
 we set $\lm f:M\to {\bf B}_E:t\mt\lm{\bf u}(t)$, where $(\lm{\bf 
u})(t):=(...,\lm u_\nu(t),...)$. Of course, for our $f$ and $\lm\in[0,1]$ we have $|||\lm 
f|||=\lm|||f|||$. 

\begin{proposition}\label{p:norm of I_Ef}
% \medskip {\bf Proposition 2.16.} {\it 
If ${\bf F}(M)$ is a $\odot$-free ${\bf L}$-space with the base $M$, then for every ${\bf L}$-space 
$E$, %  and the relevant bijection ${\bf I}_E:Set(M,\odot(E))\to\cc({\bf F}(M),E)$ 
we have $\|{\bf I}_E(f)\|_{{\bf L}b}=|||f|||$; here,  ${\bf I}_E: {\bf Set}(M,\odot(E))\to\cc({\bf F}(M),E)$ is the bijection appearing in Definition \ref{d:free space}.
\end{proposition}

\smallskip
Thus, our bijections are, in a sense, ``isometric'' maps, and to prove this we do not have to use a 
concrete construction of ${\bf F}(M)$.

\begin{proof}
% \smallskip {\it Proof}. 
First, for $\lm\in[0,1]$, applying natural property to $\va:=\lm\id_E$, and observing that 
$\odot(\lm\id_E)(f)=\lm f$, we have ${\bf I}_E(\lm f)=\lm{\bf I}_E(f)$. 

Representing every $f:M\to{\bf B}_E$ as $\lm f_0$ with $|||f_0|||=1$, we see that it suffices to 
show that the desired equality holds provided $|||f|||=1$. If it would be  not  so, then for some 
$\lm\in(0,1)$, we (still) would have $\lm^{-1}{\bf I}_E(f)\in\cc({\bf F}(M),E)$. There exists 
$g:M\to{\bf B}_E$ with ${\bf I}_E(g)=\lm^{-1}{\bf I}_E(f)$. Then we have ${\bf I}_E(\lm(g))={\bf 
I}_E(f)$. Hence $f=\lm g$ and $|||f|||=\lm|||g|||<1$, a contradiction. % $\sq$ 
\end{proof}

\begin{remark} \label{r: form of pi}
In Proposition \ref{p:quotient onto base}, we constructed the ``canonical'' morphism $\pi : {\bf F}^\sq(\sq E) \to E$. For the special case when $\sq$ is our functor $\odot$, we can describe $\pi$ in more detail.
% Recall the general-categorical ``canonical'' morphism $\pi$ from Proposition \ref{p:quotient onto base}. What is its explicit form in the case of our concrete functor $\odot$ ? 
%What is the explicit form of the morphism $\pi$ from Proposition \ref{p:quotient onto base} in the 
%case of our concrete functor $\odot$ ? 
By Proposition \ref{p:adding free}, the free ${\bf L}$-space ${\bf F}(\odot E)$ is  
 $\oplus_{\nu,u} (L^{\nu,u})^*$, where $({\nu,u})$ runs all pairs $(\nu\in\Lambda, u\in\odot E)$ and 
$(L^{\nu,u})^*$ is a copy of $(L^\nu)^*$ with the minimal ${\bf L}$-norm. 
In other words, it is the  coproduct in ${\bf LNor_1}$ of the family $\{(L^{\nu,u})^*\}$.
Here, $u$ is the collection of elements $u_\mu \in B(L^\mu \ot E)$, with $\mu$ running over $\Lambda$.
For such $(\nu,u)$, define $\va_{\nu,u}:(L^{\nu,u})^*\to E$ as ${\bf I}^\nu_E (u_\nu):(L^\nu)^*\to E$ (recall that $L^{\nu,u} = L^\nu$).
Following Proposition \ref{p:quotient onto base}, define $\pi:{\bf F}(\odot E)\to  E$ to be the coproduct $\oplus_{\nu,u}\va_{\nu,u}$,
We show directly
% ``by hand'', without applying to general-categorical  Proposition \ref{p:quotient onto base}, 
that $\pi$ is $\odot$-admissible, that is, the map $\odot\pi:\odot {\bf F}(\odot E)\to\odot E$ 
% has a right inverse $\psi$. 
is surjective.

An elements $v \in \odot {\bf F}(\odot E)$ can be represented as 
 ``rows'' $(...,v_\nu,...)$ where $v_\nu\in B(\oplus_{\mu,u}L^\nu\ot (L^{\mu,u})^*)$ (with indices $\mu\in\Lambda, u\in\odot E$); $\oplus$ refers to the $\ell_1$ sum (coproduct in ${\bf Nor_1}$, as described in the paragraph preceding Proposition \ref{p:adding free}). 
 Then $\odot \pi (v) = ( \ldots , z_\nu, \ldots)$, where
 $$
 z_\nu = \sum_{m \in \Lambda , u \in \odot E} \big( \id_{L^\nu} \otimes \oplus_{\mu,u} \phi_{\mu,u} \big) (v_\nu) \in L^\nu \otimes E .
 $$
 More specifically, we can write $v_\nu = \big(\alpha_{\mu, u}^\nu\big)_{\mu \in \Lambda, u \in \odot E}$,
with $\alpha_{\mu, u}^\nu \in L^\nu \otimes (L^{\mu,u})^*$. Then
  $$
 z_\nu = \sum_{m \in \Lambda , u \in \odot E} \big( \id_{L^\nu} \otimes \phi_{\mu,u} \big) \alpha_{\mu, u}^\nu .
 $$
 
 Now suppose we are given $z = ( \ldots , z_\nu, \ldots) \in \odot E$ (with $z_\nu \in B(L^\nu \otimes E)$ for $\nu \in \Lambda$). Our goal is to find $v \in \odot {\bf F}(\odot E)$ so that $\odot \pi (v) = z$.
 For $\nu \in \Lambda$ define by setting $v_\nu = \big(\alpha_{\mu, u}^\nu\big)_{\mu \in \Lambda, u \in \odot E}$, where $\alpha_{\nu,z_\nu}^\nu = w^\nu$ ($w^\nu$ was defined in \eqref{eq:define w}), and $\alpha_{\mu, u}^\nu = 0$ if $(\mu,u) \neq (\nu,z_n)$. Let $v = ( \ldots, v_\nu, \ldots)$.
 
 By the proof of Proposition \ref{p:norm est}, $\|w\| = 1$, hence $v \in \odot {\bf F}(\odot E)$.
 It is easy to verify that $ \big( \id_{L^\nu} \otimes \phi_{\nu,u} \big) w^\nu = u$ for any $u \in L^\nu \otimes (L^{\nu,u})^*$, hence also for $u = z_n$. Therefore, $\odot \pi(v) = z$.
\end{remark}

\section{Applications to the projectivity}\label{s:projectivity}

Suppose, for a moment, that we have an arbitrary base space and an arbitrary functor $\sq:{\bf 
L}{\bf Nor_1}\to{\bf Set}$. 

\begin{definition}\label{d:generalized projectivity}
% {\bf Definition 3.1.} 
An ${\bf L}$-space $P$ is called {\it $\sq$-projective}, if for every ${\bf L}$--spaces $G$ and 
$E$, every $\sq$-admissible operator $\tau:G\to E$ and every  ${\bf L}$-contractive operator 
$\va:P\to E$ there exists an ${\bf L}$-contractive lifting of $\va$ across $\tau$. Our $P$ is 
called {\it asymptotically $\sq$--projective}, if for every $G,E$ and $\tau$ as before, and 
$\va:P\to E$ with $\|\va\|_{{\bf L}b}<1$ there exists an ${\bf L}$--contractive lifting $\psi$ of 
$\va$ across $\tau$. 
\end{definition}

\begin{proposition}\label{p:retract is projective}
% \medskip {\bf Proposition 3.2.} {\it 
Every retract of a $\sq$-projective ${\bf L}$-space is itself a  $\sq$-projective ${\bf L}$-space. Every near-retract of an asymptotically $\sq$-projective ${\bf L}$--space is itself an asymptotically $\sq$-projective ${\bf L}$--space. % } 
\end{proposition}

\begin{proof}
% \smallskip {\it Proof.}
% The proof repeats the proof of
Proceed as in Proposition \ref{p:retract of projective}, with obvious modifications. % $\sq$
\end{proof}

\begin{proposition}\label{p:free spaces are projective}
% \medskip {\bf Proposition 3.3.} {\it 
$(i)$ For every $M$, the $\sq$-free ${\bf L}$-space ${\bf F}^\sq(M)$ $($if it does exist$)$ is 
$\sq$-projective. 

\smallskip
$(ii)$ A retract, respectively near-retract of a $\sq$-free ${\bf L}$-space is $\sq$--projective, 
respectively asymptotically $\sq$--projective. 

\smallskip
$(iii)$ If every set is a base of some $\sq$-free ${\bf L}$-space, then every $\sq$--projective ${\bf 
L}$-space is a retract, and every asymptotically $\sq$--projective ${\bf L}$-space a near-retract 
of some $\sq$-free ${\bf L}$-space. 
% }
%there exists a $\sq$-free space, say ${\bf F}$, with the base $\sq E$, than every 
%$\sq$--projective ${\bf L}$-space is a retract, and every asymptotically $\sq$--projective ${\bf 
%L}$-space a near-retract of ${\bf F}$.} 
\end{proposition}

% \smallskip {\it Proof}. 
\begin{proof}
(i) It is a particular case of a known categorical statement. 
%Nevertheless we give an outline of the proof. 
Suppose $\tau$ is admissible; find a map $\rho:\sq E\to\sq G$ so that $(\sq\tau)\rho=\id_{\sq E}$.
For $\va:{\bf F}^\sq(M)\to E$ let $\psi= {\bf I}_G(\rho \circ {\bf I}_E^{-1}(\va))$, then $\tau\psi=\va$.
 %\query{I have switched ${\bf I}^\sq$ to ${\bf I}$ for the sake of notational uniformity.}

\smallskip
(ii) follows from p.~(i), Proposition \ref{p:retract is projective}, and the obvious observation that $\sq$--projective ${\bf L}$-spaces are asymptotically $\sq$--projective. 

\smallskip
(iii) Let $P$ be a $\sq$-projective ${\bf L}$-space. Proposition \ref{p:quotient onto base} provides a $\sq$--admissible operator $\pi: {\bf F}(\sq P) \to P$. By the projectivity of $P$, $\va = \id_P$ has an ${\bf L}$-contractive lifting across $\pi$, which we denote by $\psi$. In other words, $\psi$ is an ${\bf L}$-contractive right inverse of $\pi$ -- that is, $P$ is a retract of ${\bf F}(\sq P)$.
% Taking $\va:=\id_P$, we see that $\pi$ is a retraction.

If $P$ is only asymptotically projective, the preceding proof needs to be modified.
For any $\e \in (0,1)$, $\va:= (1+\e)^{-1} \id_P$ has an ${\bf L}$-contractive lifting $\psi$ across $\pi$. Then $(1+\e)\psi$ is a right inverse of $\pi$, hence $\pi$ is a near retraction.
\end{proof}

% \medskip
Return to our special functors $\odot,\odot^0:{\bf L}{\bf Nor_1}\to{\bf Set}$, corresponding to a 
(so far arbitrary) family of contractively complemented subspaces $L^\nu$ of ${\bf L}$. Propositions 
\ref{p:coisometric admissible} and \ref{p:odot criterion}  imply:

\begin{proposition}\label{p:odot projectivity}
In the above notation, any $\odot$-projective $(\odot^0$-projective$)$ ${\bf L}$-space is metrically projective $($respectively, extremely projective$)$.
If ${\bf L}$ is properly presented relative to the family $L^\nu : \nl$, then the converse impications also hold.
% and the functors $\odot, \odot^0$ are constructed based on this family. Then an ${\bf L}$-space is 
\end{proposition}

% From now on we assume that the spaces $L^\nu$ are finite dimensional, and ${\bf L}$ is properly presented.
% This is the case, for instance, when ${\bf L} = L_p^0(X)$, and the family $(L^\nu)$ is as described in in Example \ref{e:contractively complemented subspaces}.

% {\bf Remark 3.6.} 
Next we establish a connection between projectivity and freeness. Part (ii) of Theorem \ref{t:F(M)} below can be deduced from the general theory of asymptotic rigged categories, introduced in~\cite{hemsb, heruv}. For the sake of clarity, we present a direct approach.  
% Actually, p. (ii) can be deduced from a certain general-categorical assertion, concerning so-called asymptotic rigged categories, introduced in~\cite{hemsb, heruv}. But here, for the sake of clarity, we prefer a direct approach.  

\begin{theorem}\label{t:F(M)}
% \medskip {\bf Theorem 3.5.} {\it 
Let $M$ be an arbitrary set, and  ${\bf F}(M)$ the $\odot$--free ${\bf L}$-space with the base $M$ 
$($cf. Theorem \ref{t:criterion for freeness}$)$. Then:

\smallskip
(i) ${\bf F}(M)$ is $\odot$-projective and metrically projective. 

\smallskip
(ii) ${\bf F}(M)$ is asymptotically $\odot^0$-projective and  extremely projective. 
% }
\end{theorem}
 
 \begin{proof}
% \smallskip {\it Proof}. 
$(i)$ follows from Propositions \ref{p:free spaces are projective}(i) and \ref{p:odot 
projectivity} combined.   
 
 \smallskip
$(ii)$  Let $G$ and $E$ be ${\bf L}$-spaces, $\va:{\bf F}(M)\to E$ a completely contractive operator 
with $\|\va\|_{{\bf L}b}=:\theta<1$, and $\tau:G\to E$ an $\odot^0$-admissible operator. 
  According to Propositions \ref{p:characterization of projectivity} and \ref{p:odot projectivity}, our task is to find an ${\bf L}$-contractive operator $\psi:{\bf F}(M)\to G$,   making the diagram \eqref{eq:lifting} commutative (that is, $\tau \psi = \phi$). 

By the condition on $\tau$, for every $\nl$ the restriction of $\tau_\ii$ to $L^\nu\ot G$ is a 
coisometric operator. Define the map $\rho_\nu : B^\nu_E \to B^\nu_G$: if $u \in B^\nu_E$ satisfies $\|u\| \leq \theta$, find $\rho_\nu(u)$ so that $\tau_\ii \rho_\nu(u) = u$. If $\|u\| > \theta$, let $\rho_\nu(u)$ be an arbitrary element of $B^\nu_G$.
% Therefore for  $u\in B^\nu_E$ we can fix  an arbitrary $v\in B^\nu_G$ such 
% that $\tau_\ii(v)=u$ provided $\|u\|\le\theta$ and (just) an arbitrary $v$ otherwise. In both cases we denote this $v$ by $\rho_\nu(u)$ and
Now consider the map $\rho:{\bf B}_E\to{\bf B}_G$, taking ${\bf u}=(...,u_\nu,...)$ to ${\bf v}=(...,\rho_\nu(u_\nu),...)$. We see, in the notation of Proposition \ref{p:norm of I_Ef}, that if ${\bf u}:M\to{\bf B}_E$ 
satisfies $|||{\bf u}|||\le\theta$, then 
\[
(\odot\tau)\rho({\bf u}(t))={\bf u}(t). %\eqno(5) 
\] 
Now replace ${\bf u}$ with ${\bf I}_E^{-1}(\va):M\to{\bf B}_E$, where ${\bf I}_E:{\bf Set}(M,{\bf B}_E)\to\cc({\bf F}(M),E)$ is the bijection from Definition \ref{d:free space}. Taking into account the aforementioned proposition and \eqref{eq:oplus norms}, for every $t\in M$ we obtain 
\[
(\odot\tau)\rho[{\bf I}^{-1}_E(\va)]= {\bf I}_E^{-1}(\va). 
\]

Finally, set $\psi:={\bf I}_G(\rho[{\bf I}_E^{-1}(\va)])$. As our bijections are natural, we conclude that $\tau\psi={\bf I}_E[(\odot\tau)\rho[{\bf I}_E^{-1}(\va)]]=\va$. 
% With the help of Proposition \ref{p:odot projectivity}, we are done. % $\sq$ 
\end{proof}

% Recall that ${\bf L}$ is properly presented, if it is the set-theoretical union  of a family of its contractively complemented finite-dimensional subspaces.
Suppose now that $L^\nu ; \nl$ is the family of all contractively complemented finite dimensional subspaces of ${\bf L}$, and ${\bf L}$ is properly presented relative to $(L^\nu)$
(this is the case, for instance, when ${\bf L} = L_p^0(X)$, and the family $(L^\nu)$ is as described in in Example \ref{e:contractively complemented subspaces}).
Let $\odot$ be the functor corresponding to $L^\nu ; \nl$.
 Proposition \ref{p:free spaces are projective}(i,iii), considered for the case $\sq:=\odot$, together  with Proposition \ref{p:odot projectivity} gives

 \begin{proposition}\label{p:conclusion about properly presented}
% \medskip {\bf Proposition 3.7.} {\it 
% Suppose that ${\bf L}$ is properly presented (like in the case ${\bf L}:=L_p^0(X)$), and the functor $\odot$ corresponds to the family of all relevant subspaces. Then
In the above notation, every metrically projective ${\bf L}$--space is a retract, and every extremely projective ${\bf L}$-space a near-retract of some $\odot$--free ${\bf L}$--space. $\sq$ 
\end{proposition}

\medskip
Now, combining this proposition with the explicit construction of the $\odot$-free ${\bf L}$--space 
${\bf F}(M)$, provided by Propositions \ref{p:singleton base} and \ref{p:adding free}, we obtain 
the following description of projective ${\bf L}$-spaces. 

% Let us call an ${\bf L}$-space {\it well composed}, if it is the $\oplus_1$-sum of some family of ${\bf L}$-spaces such that each of them is minimal ${\bf L}$-space $(L)^*$ for some finite-dimensional contractively complemented subspace $L$ of ${\bf L}$. (These $L$ can be different for different summands.)

We say that an ${\bf L}$-space is {\it well composed} if it is the $\ell_1$ sum 
$\oplus_{\mu \in \mathcal M} {\rm{MIN}}(Z_\mu^*)$. Here $Z_\mu$ is a finite dimensional contractively complemented subspace of ${\bf L}$, and ${\rm{MIN}}(Z_\mu^*)$ refers to the minimal ${\bf 
L}$ structure of its dual -- that is, we equip ${\bf L} \ot Z_\mu^*$ with the injective tensor norm.
Proposition \ref{p:adding free} shows that $\odot$-free spaces are well composed; moreover, for every $\mu$ there exists $\nl$ so that $Z_\mu = L^\nu$, and the cardinality of $\{\mu \in \mathcal M : Z_\mu = L^\mu\}$ is independent of $\nu$.

% Note that, by the construction of $\odot$-free spaces, given in Proposition \ref{p:adding free}, shows that such a space is $\odot$-free exactly when the cardinality of the set of summands $L_\nu^*$ is the same for every $\nu\in\Lambda$.)
%${\bf F}^\nu(\star)$ is the same for every $\nu$. 

\begin{theorem}\label{t:describe_projectivity}
% \medskip {\bf Theorem 3.8.} {\it 
$(i)$ For every ${\bf L}$ a  retract, respectively near-retract of a well composed ${\bf L}$--space is metrically, respectively extremely projective. 

$(ii)$ If ${\bf L}$ is properly presented, then every metrically, respectively extremely projective ${\bf L}$--space is a retract, respectively near-retract of some well composed ${\bf L}$--space. 
% } 

$(iii)$ If ${\bf L}$ is properly presented, then every ${\bf L}$--space is a strictly coisometric image of a well-composed ${\bf L}$--space.
\end{theorem}

\begin{proof}
% \smallskip {\it Proof.} 
(i) We provide a proof for metric projectivity, as extreme projectivity is handled similarly.
Suppose $\tau : G \to E$ is a strict ${\bf L}$-coisometry. It suffices to show that, if $L$ is a finite dimensional subspace of ${\bf L}$, complemented via a contractive projection $Q$, then any ${\bf L}$-contractive operator $\va : L^* \to E$ admits an ${\bf L}$-contractive lifting $\psi : L^* \to G$, with $\tau\psi = \va$.

For an ${\bf L}$-space $U$, define the map ${\bf I}_U : {\bf L}\ot U\to {\mathcal L}(L^*,U)$ in a manner similar to ${\bf I}_U^\nu$, introduced before Proposition \ref{p:norm est}).
For $u\in{\bf L}\ot U$, ${\bf I}_U u$ is the operator $L^* \ni g \mapsto (g\ot\id_U)(Q\cd u)$. 
Proposition \ref{p:norm est} shows that $\|u\| \geq \| {\bf I}_U u\|_{{\bf L}b}$ fr any $u \in {\bf L} \otimes E$; moreover, ${\bf I}_U$ implements a bijective isometry from $L \otimes U$ (regarded as a subspace of ${\bf L} \otimes U$) onto $\cb(L^*,U)$.
Abusing the notation slightly, we talk about ${\bf I}_U^{-1} : \cb(L^*,U) \to L \otimes U$.

Now fix $\va \in \cc(L^*,E)$. By the coisometric property of $\tau_\infty$, there exists $u \in B({\bf L} \otimes G)$ so that $\tau_\infty u = {\bf I}_E^{-1} \va$. Clearly $Q \cdot u$ has the same properties. Then $\psi := {\bf I}_G (Q \cdot u)$ is the desired lifting.

% follows from the obvious observation that every well composed ${\bf L}$-space is retract of some free ${\bf L}$-space, together with Theorem \ref{t:F(M)} and Proposition \ref{p:retract of projective}. 

(ii) is a direct corollary of Proposition \ref{p:conclusion about properly presented}. 

(iii) Suppose $E$ is an ${\bf L}$-space. Use Theorem \ref{t:criterion for freeness} to find the free object ${\bf F}(\odot E)$. By Proposition \ref{p:quotient onto base}, there exists an $\odot$-admissible map from ${\bf F}(\odot E)$ onto $E$. By Propositions \ref{p:odot criterion}, such a map is ${\bf L}$-strongly coisometric.
% follows by combining Propositions \ref{p:odot criterion} and \ref{p:quotient onto base} with Theorem \ref{t:criterion for freeness}.
\end{proof}

Next we specialize to the case when ${\bf L}$ is either $L_p(X)$ or $L_p^0(X)$. In the discussion preceding Definition \ref{d:L space}, we established that the classes of $L_p(X)$-spaces and $L_p^0(X)$-spaces are identical (that is, an $L_p(X)$-space must be an $L_p^0(X)$-space, and vice versa; there is a bijective correspondence between norms as well). 
% We shall say that an $L_p^0(X)$-space is {\it nice} if it is an $\ell_1$ sum $\oplus_{\mu \in \mathcal M} {\rm{MIN}}\big(\ell_q^{n_\mu}\big)$, for some family $(n_\mu)_{\mu \in {\mathcal{M}}}$; here, $1/p + 1/q = 1$.
We begin our analysis by stating a lemma, which links the notions of ${\bf L}$-(strict) coisometricity, and being a (near) retract, for bases ${\bf L} = L_p(X)$ and ${\bf L} = L_p^0(X)$.

\begin{lemma}\label{l:LpX vs Lp0X projectvity}
Suppose $E$ and $G$ are $L_p^0(X)$-spaces, or, equivalently, $L_p(X)$-spaces.

$(i)$ If $\tau : G \to E$ is an $L_p(X)$-coisometric $(L_p(X)$-strictly coisometric$)$ map, then it is also $L_p^0(X)$-coisometric $($respectively, $L_p^0(X)$-strictly coisometric$)$.

$(ii)$ If $\tau : G \to E$ is $L_p^0(X)$-coisometric, then it is also $L_p(X)$-coisometric. 

$(iii)$ A map $\tau : G \to E$ is a retraction (near-retraction) in the category of $L_p(X)$-spaces if and only if it is a retraction (near-retraction) in the category of $L_p^0(X)$-spaces.
\end{lemma}

Remark \ref{r:coisometries} below shows that part (ii) of the lemma cannot be improved. Specifically, it provides an example of an $L_p^0(X)$-strictly coisometric map $\tau : G \to E$, which is not $L_p(X)$-strictly coisometric.

\begin{proof}
(i) We deal with strict coisometries (the case of coisometries is similar). Fix $u = \sum_{k=1}^N \xi_k \otimes u_k \in L_p^0(X) \otimes E \subset L_p(X) \otimes E$. Due to $\tau$ being $L_p(X)$-strictly coisometric, there exists $v \in L_p(X) \otimes G$ so that $\tau_\infty v = u$, and
$$\|v\|_{L_p(X) \otimes G} = \|u\|_{L_p(X) \otimes E} = \|u\|_{L_p^0(X) \otimes E} . $$
Find $\nl$ so that $\xi_k \in L^\nu$ for all $k$, and let $v' = Q^\nu v$. Clearly $v' \in L^\nu \otimes G \subset L_p^0(X) \otimes G$, $\tau_\infty v' = u$, and
$$ \|v'\|_{L_p^0(X) \otimes G} = \|v'\|_{L_p(X) \otimes G} \leq \|v\| = \|u\| , $$
establishing the $L_p^0(X)$-strict coisometricity of $\tau$.

(ii) Consider $u = \sum_{k=1}^N \xi_k \otimes x_k \in L_p(X) \otimes E$. Fix $\varepsilon > 0$; our goal is to find $v \in L_p(X) \otimes G$ so that $\|v\| < \|u\| + \varepsilon$, and $\tau_\infty v = u$. To this end, find $\nl$ so that $\sum_k \|\xi_k - Q^\nu \xi_k\| \|x_k\| < \varepsilon/2$. Further, find $w \in L_p^0(X) \otimes G$ so that $\tau_\infty w = \sum_k Q^\nu \xi_k \otimes x_k$, and $\|w\| < \|u\| + \varepsilon/2$. For each $k$ find $y_k \in G$ so that $\tau y_k = x_k$, and $\sum_k \|\xi_k - Q^\nu \xi_k\| \|y_k\| < \varepsilon/2$. Then $v = w + \sum_k \big( \xi_k - Q^\nu \xi_k \big) \otimes y_k$ is the desired lifting of $u$.

(iii) This statement follows from the fact that, for any map $\va$ between $L_p(X)$-spaces (or, equivalently, $L_p^0(X)$-spaces), $\|\va\|_{L_p(X) b} = \|\va\|_{L_p^0(X) b}$.
\end{proof}

Note that $L_p^0(X)$ is properly presented, and the spaces $L^\nu$ involved in this proper presentation are isometric copies of $\ell_p^{n_\nu}$. This immediately leads to:

\begin{lemma}\label{l:well composed}
An $L_p^0(X)$-space is well-composed if and only if it is an $\ell_1$ sum $\oplus_{\mu \in \mathcal
M} {\rm{MIN}}\big(\ell_q^{n_\mu}\big)$, for some family $(n_\mu)_{\mu \in {\mathcal{M}}}$; here, $1/p + 1/q = 1$.   $\sq$ 
\end{lemma}

Having laid the groundwork, we can now characterize projective $L_p(X)$ and $L_p^0(X)$ spaces.

\begin{theorem}\label{t:projective LpX}
% For $1 \leq p \leq \infty$, find $q$ so that $1/p + 1/q = 1$.
$(i)$ An $L_p^0(X)$-space is metrically $($extremely$)$ projective if and only if it is a retract $($respectively, near retract$)$ of a well-composed $L_p^0(X)$-space.
% an $\ell_1$ sum $\oplus_{\mu \in \mathcal M} {\rm{MIN}}\big(\ell_q^{n_\mu}\big)$, for some family $(n_\mu)_{\mu \in {\mathcal{M}}}$.

$(ii)$ An $L_p(X)$-space is extremely projective if and only if it is a near retract of a well-composed $L_p^0(X)$-space.
% an $\ell_1$ sum $\oplus_{\mu \in \mathcal M} {\rm{MIN}}\big(\ell_q^{n_\mu}\big)$, for some family $(n_\mu)_{\mu \in {\mathcal{M}}}$.

$(iii)$ Any $L_p(X)$-space which is a retract of a well-composed $L_p^0(X)$-space
% an $\ell_1$ sum $\oplus_{\mu \in \mathcal M} {\rm{MIN}}\big(\ell_q^{n_\mu}\big)$
is metrically projective.

$(iv)$ Any $L_p^0(X)$-space $(L_p(X)$-space$)$ is an image of a well-composed $L_p^0(X)$-space under an $L_p^0(X)$-strict coisometry $($respectively, $L_p(X)$-coisometry$)$.
%
%  Suppose ${\bf L}$ is either $L_p(X)$ or $L_p^0(X)$. Then an ${\bf L}$-space $E$ is metrically $($extremely$)$ projective if and only if it is a retract $($respectively, near retract$)$ of an $\ell_1$ sum $\oplus_{\mu \in \mathcal M} {\rm{MIN}}\big(\ell_q^{n_\mu}\big)$, for some family $(n_\mu)_{\mu \in {\mathcal{M}}}$. Here, $1/p + 1/q = 1$.
\end{theorem}

Proposition \ref{p:not a strict coisometry} below shows that, in part (iv), an $L_p(X)$-space need not be a strictly coisometric image of a well-composed space.

% We begin our proof of Theorem \ref{t:projective LpX} by establishing Lemma \ref{l:LpX vs Lp0X projectvity}, which allows us to connect ${\bf L}$-(strict) coisometricity for bases ${\bf L} = L_p(X)$ and ${\bf L} = L_p^0(X)$.

\begin{proof}%[Proof of Theorem \ref{t:projective LpX}]
(i) follows from Theorem \ref{t:describe_projectivity}(i,ii).
% Therefore, it suffices to consider the case of ${\bf L} = L_p^0(X)$.
% The spaces $L^\nu$ involved in the proper presentation are isometric copies of $\ell_p^{n_\nu}$. It remains to apply Theorem \ref{t:describe_projectivity}(i,ii). 

(ii) By Lemma \ref{l:LpX vs Lp0X projectvity}(i,ii), $P$ is extremely projective as an $L_p(X)$-space if and only if it is extremely projective as an $L_p^0(X)$-space. By part (i) of this theorem, the last condition is equivalent to being a near-retract of a well-composed $L_p^0(X)$-space. By Lemma \ref{l:LpX vs Lp0X projectvity}(iii), the notions of being a retract or near-retract in the categories of $L_p(X)$-spaces and $L_p^0(X)$-spaces coincide.

% and (iii) follow by combining (i) with Lemma \ref{l:LpX vs Lp0X projectvity}.

(iii) Suppose $E, G, P$ are $L_p(X)$-spaces, $P$ is a retract of a well-composed $L_p^0(X)$-space, $\va : P \to E$ is an $L_p(X)$-contraction, and $\tau : G \to E$ is an $L_p(X)$-strict coisometry. By Lemma \ref{l:LpX vs Lp0X projectvity}(i), $\tau$ is an $L_p^0(X)$-strict coisometry as well. Combining part (i) of this theorem with Lemma \ref{l:LpX vs Lp0X projectvity}(iii), we conclude that $\va$ has an $L_p^0(X)$-contractive lifting $\psi : P \to G$. To complete the proof, note that $\psi$ must be $L_p(X)$-contractive as well.

In (iv), the first statement is a consequence of Theorem \ref{t:describe_projectivity}(iii). The second one follows by observing that any $L_p^0(X)$-strict coisometry is an $L_p(X)$-coisometry.
\end{proof}

\begin{remark}\label{r:coisometries}
We do not know whether a metrically projective $L_p(X)$-space is necessarily metrically projective in the category of $L_p^0(X)$-spaces.
The following example shows that an $L_p^0(X)$-strict coisometry may fail to be an $L_p(X)$-strict coisometry.

Let $I = [0,2\pi)$. Equip $G = \ell_1(I)$ and $E = \ell_2^2$ with an $L_1$-structure arising from $L_1(G)$ and $L_1(E)$ respectively (we use $L_1$ and $L_1^0$ for $L_1(0,1)$ and $L_1^0(0,1)$). These are indeed $L_1$-spaces, by \cite[Sections 1.8-9]{dal}. Denote the canonical bases of $G$ and $E$ by $g_t$ ($t \in I$) and $e_1, e_2$, respectively.

Define $\tau : G \to E : g_t \mapsto \cos t e_1 + \sin t e_2$, and note it is a strict coisometry on the Banach space level. Indeed, any $x \in E$ can be written as $\|x\| (\cos t e_1 + \sin t e_2)$, with some $t \in I$. Then $x = \tau( \|x\| g_t)$.

Next show that $\tau$ is $L_1^0$-strictly coisometric. Indeed, any element $u \in L_1^0 \otimes E$ can be written as $u = \sum_k \id_{A_k} \otimes x_k$, with $x_k \in E$, and $A_k$ being disjoint measurable subsets of $(0,1)$. Then $\|u\| = \sum_k \lambda(A_k) \|x_k\|$, where $\lambda$ is the Lebesgue measure. For each $k$ find $y_k \in G$ with $\|x_k\| = \|y_k\|$, $\tau y_k = x_k$. Let $v = \sum_k \id_{A_k} \otimes x_k$, and note that $\tau_\infty v = u$, and $\|v\| = \|u\|$.

To show that $\tau$ is not $L_1$-strictly coisometric, note first that, if $\|y\| = 1$, and $\tau y = \cos t e_1 + \sin t e_2$ (for some $t \in I$), then $y = g_t$. Indeed, write $y = \sum_{s \in I} \alpha_s g_s$, where $\sum_s |\alpha_s| = 1$ (hence $\alpha_s = 0$ for at most countably many values of $s$). Then
$$
1 = \langle \cos t e_1 + \sin t e_2, \tau y \rangle + \sum_s \alpha s \langle \cos t e_1 + \sin t e_2 , \cos s e_1 + \sin s e_2 \rangle = \sum_s \alpha_s \cos (t-s) . 
$$
The only way for this equality to hold is to have $\alpha_t = 1$, $\alpha_s = 0$ for $s \neq t$.

Now consider $u = \xi_1 \otimes e_1 + \xi_2 \otimes e_2$, where $\xi_1(s) = \cos 2 \pi s$ and $\xi_2(s) = \sin 2 \pi s$ for $s \in (0,1)$. Thus, $u(s) = \cos (2 \pi s) e_1 + \sin (2 \pi s) e_2$ for any $s$. Consequently, $\|u\|_{L_1(E)} = \int_0^1 \| u(s)\| \, ds = 1$. We shall show that $u$ has no norm $1$ lifting.

Suppose, for the sake of contradiction, that there exists $v = \sum_{k=1}^N \xi_k \otimes y_k \in L_1 \otimes G$ so that $\tau_\infty v = u$, and $\|v\| = 1$. In the coordinate form, we must have that $\| \sum_k \xi_k(s) y_k \| = 1$, and $\tau \big( \sum_k \xi_k(s) y_k \big) = u(s) = \cos (2 \pi s) e_1 + \sin (2 \pi s) e_2$, for any $s \in {\mathcal{S}}$, where ${\mathcal{S}}$ has measure $1$. Consequently, $\sum_k \xi_k(s) y_k = g_{2 \pi s}$, for any $s \in {\mathcal{S}}$. This, however, is impossible, as each $y_k$ has countable support in $I$.
\end{remark}

We close this paper by an example indicating that the category of $L_p(X)$-spaces may not be the ``right'' one.

\begin{proposition}\label{p:not a strict coisometry}
 There exists an $L_1(0,1)$-space which cannot be represented as an $L_1(0,1)$-strictly coisometric image of a well-composed $L_1^0(0,1)$-space.
\end{proposition}

\begin{proof}
 Denote by $E$ the space $C(\T)$, where $\T$ is the unit circle, equipped with its minimal $L_1(0,1)$ structure -- that is, the norm on $L_1(0,1) \otimes C(\T) \simeq C(\T, L_1(0,1))$ comes from the injective tensor product.
%  Suppose, for the sake of contradiction, that there is a strict coisometry $\tau$ from the $\ell_1$ sum $\oplus {\rm{MIN}}(\ell_1^{n_i})$; we represent $t = \oplus \tau_i$, with $\tau_i$ acting on the summand $\ell_1^{n_i}$.
 Consider normalized independent Gaussian random variables $g_1, g_2 \in L_1(0,1)$, and the coordinate functions $f_1, f_2 \in C(\T)$. Let $u = g_1 \otimes f_1 + g_2 \otimes f_2$. Then $\|u\| = 1$. In fact, $u$ corresponds to the weak$^*$ to norm continuous operator
 $$
 \tilde{u} : L_\infty(0,1) \to C(\T) : h \mapsto \Big( \int h g_1 \Big) f_1 + \Big( \int h g_2 \Big) f_2 .
 $$
 One can write $\tilde{u} = j q$; here $j : \ell_2^2 \to C(\T) : \delta_i \mapsto f_i$ ($\delta_1, \delta_2$ form the canonical basis in $\ell_2^2$) is an isometric embedding, and $q$ is the strictly coisometric adjoint of $q_* : \ell_2^2 \to L_1 : \delta_i \mapsto g_i$. Consequently, $\tilde{u}$ maps the closed unit ball of $L_\infty$ onto the closed unit ball of $j(\ell_2^2) \subset L_1$.
 
  Suppose, for the sake of contradiction, that there exists an $L_1(0,1)$-strict coisometry $\tau$ from an $\ell_1$ sum $\oplus {\rm{MIN}}(\ell_\infty^{n_i})$ onto $C(\T)$. Find $v \in L_1(0,1) \otimes \oplus {\rm{MIN}}(\ell_\infty^{n_i})$ so that $\|v\|=1$, and $\id_{L_1} \otimes \tau(v) = u$. 
 As $\oplus$ denotes the algebraic sum, we can write $v$ as a finite sum $\oplus_{i \in I} v_i$, with $v_i \in L_1(0,1) \otimes \ell_\infty^{n_i}$, and $\|v\| = \sum_i \|v_i\|$. 
%  To each $v_i$ corresponds an operator $\tilde{v}_i : L_\infty(0,1) \to \ell_1^{n_i}$. Also, write $\tau = \oplus \tau_i$, with $\tau_i$ being the restriction of $\tau$ to $\ell_1^{n_i}$. Clearly $\|\tau_i\| \leq \|\tau\| \leq 1$. 
Based on $v$, we construct an operator $\tilde{v}$ from $L_\infty(0,1)$ into the finite dimensional space $Z = \big( \oplus_{i \in I} \ell_\infty^{n_i} \big)_1$. The equality $u = (\id_{L_1} \otimes \tau) v$ translates into $\tilde{u} = \tau|_Z \tilde{v}$. 

As noted above, for every $\eta \in \ran \tilde{u}$ we can find $\xi \in L_\infty(0,1)$ so that $\|\xi\| = \|\eta\|$, and $\tilde{u} \xi = \eta$. Now let $Y = \ran \tilde{v} \subset Z$; then, $\tau(Y) \subset \ran \tilde{u}$, and for every $\eta \in \ran \tilde{u}$ we can find $y \in Y$ so that $\|y\| = \|\eta\|$, and $\tau y = \eta$. In other words, $\tau$ acts as a strict coisometry from $Y$ onto $\ran \tilde{u} \simeq \ell_2^2$.
 However, the unit ball of $Z$ is a polytope, hence the same is true for $Y$. As a circle is not a cross-section of a polytope, no coisometry from $Y$ onto a copy of $\ell_2^2$ exists.
\end{proof}

\end{document}